\documentclass[final,3p]{elsarticle}

\usepackage[colorlinks=true]{hyperref}
\usepackage[hyperref,dvipsnames]{xcolor}
\usepackage[normalem]{ulem}
\usepackage{cancel}
\usepackage[T1]{fontenc} 							
\usepackage[english]{babel} 					
\usepackage{amssymb,amsmath,amsthm,amsbsy}		
\usepackage{array} 						
\usepackage{pdfpages}								
\usepackage{graphicx} 					
\usepackage{caption}[2004/07/16]					
\usepackage{subfigure}							
\usepackage{mathrsfs} 	
\usepackage{empheq}								
\usepackage{nicefrac}								
\usepackage{stmaryrd}								
\usepackage{fixltx2e} 							
\usepackage{fix-cm}								
\usepackage{bm}
\usepackage{booktabs}
\usepackage{here}								
\usepackage{url}
\usepackage{enumitem} 
\usepackage{thmbox}
\usepackage{shadethm}		
\usepackage{siunitx}  		
\usepackage{titlesec}	 

\usepackage{etoolbox}
\robustify\bfseries
\newcommand{\uu}{{\boldsymbol{u}}}
\newcommand{\vv}{{\boldsymbol{v}}}
\newcommand{\ww}{{\boldsymbol{w}}}
\newcommand{\ff}{{\boldsymbol{f}}}

\newcommand{\x}{{\boldsymbol{x}}}
\newcommand{\g}{{\boldsymbol{g}}}
\newcommand{\n}{{\boldsymbol{n}}}
\newcommand{\s}{{\boldsymbol{s}}}
\newcommand{\drm}{{\mathrm{d}}}

\newcommand{\Pra}{{\mathrm{Pr}}}										
\newcommand{\Ray}{{\mathrm{Ra}}}	
\newcommand{\Gra}{{\mathrm{Gr}}}	
\newcommand{\Nu}{{\mathrm{Nu}}}	
\newcommand{\NuAvg}{{\overline{\mathrm{Nu}}}}								
\newcommand{\hot}{{\mathrm{hot}}}							
\newcommand{\cold}{{\mathrm{cold}}}
\newcommand{\refrm}{{\mathrm{ref}}}
\newcommand{\midrm}{{\mathrm{mid}}}
\newcommand{\maxrm}{{\mathrm{max}}}
\newcommand{\kpr}{{k^\prime}}
\newcommand{\dvg}{{\mathrm{div}}}
\newcommand{\RT}{{\mathcal{RT}}}

\newcommand\rb[1]{\hspace{-0.1ex}\left(#1\right)}
\newcommand\sqb[1]{\left[#1\right]}
\newcommand\abs[1]{\left\lvert#1\right\rvert}
\newcommand\norm[1]{\left\lVert#1\right\rVert}
\newcommand{\otoprule}{\midrule[\heavyrulewidth]}
\newcommand\Area[1]{\mathcal{A}\rb{#1}}
\newcommand\Vol[1]{\mathcal{V}\rb{#1}}
	
\newcommand{\Q}{Q_{\mathrm{IF}}}
\newcommand{\tend}{t_{\mathrm{end}}}
\newcommand{\sip}{{\mathrm{sip}}}
\newcommand{\upw}{{\mathrm{upw}}}
\newcommand{\eps}{{\varepsilon}}
\newcommand{\R}{\mathbb{R}}
\newcommand{\tr}{{\mathrm{tr},k}}

\newcommand{\T}{{\mathcal{T}_h}} 									
\newcommand{\F}{{\mathcal{F}_h}}										
\newcommand{\Fi}{{\mathcal{F}_h^i}}									
\newcommand{\Fb}{{\mathcal{F}_h^b}}									
\newcommand{\FbD}{{\mathcal{F}_h^{D}}}								
\newcommand{\FbN}{{\mathcal{F}_h^{N}}}
\newcommand{\FK}{{\mathcal{F}_K}}
								
\newcommand{\lavg}{\big\{\hspace{-0.95ex}\big\{}						
\newcommand{\ravg}{\big\}\hspace{-0.99ex}\big\}}						
\newcommand{\ljmp}{\left\llbracket}									
\newcommand{\rjmp}{\right\rrbracket}									
\newcommand\jmp[1]{\ljmp#1\rjmp}										
\newcommand\avg[1]{\lavg#1\ravg}		

\newcommand\restr[2]{{												
	\left.\kern-\nulldelimiterspace									
	#1
	\vphantom{\big|}
	\right|_{#2}
	}}
	
\newcommand\RTquad[1]{{												
	\mathcal{RT}_{\sqb{#1}}
	}}
	
\newcommand\Quad[1]{{												
	\mathcal{Q}_{#1}
	}}

\newcommand{\goodgap}{%
	\hspace{\subfigtopskip}
	\hspace{\subfigbottomskip}
	}							

\newcommand{\Lp}{{\mathcal{L}}}		
\newcommand{\Hk}{{\mathcal{H}}}		
\newcommand{\Hdiv}{\boldsymbol{\mathcal{H}}\rb{\dvg;\Omega}}								
\newcommand{\V}{{\mathcal{V}}}										
\newcommand{\UU}{{\boldsymbol{\mathcal{U}}}}							
\newcommand{\PP}{{\mathcal{P}}}

\titleformat{\section}{\fontsize{10.5}{17}\bfseries}{\thesection}{1em}{}
\titleformat{\subsection}{\bfseries\itshape}{\thesubsection}{1em}{}		
		
\journal{the `Journal of Computational Physics' (accepted: January 25, 2017)}

\definecolor{mediumblue}{RGB}{0,0,205}
\definecolor{forestgreen}{RGB}{34,139,34}
\definecolor{darkred}{RGB}{200,0,0}

\begin{document}

\hypersetup{
  linkcolor=darkred,
  urlcolor=forestgreen,
  citecolor=mediumblue
}

\begin{frontmatter}


\title{Stabilised dG-FEM for incompressible natural convection flows with boundary and moving interior layers on non-adapted meshes}

\author[]{Philipp W.\ Schroeder\corref{cor1} \fnref{fn1}}
	\cortext[cor1]{Corresponding author}
	\fntext[fn1]{ORCID: \url{https://orcid.org/0000-0001-7644-4693}}	
    \ead{p.schroeder@math.uni-goettingen.de}

\author[]{Gert Lube} 
	\ead{lube@math.uni-goettingen.de}

\address{Institute for Numerical and Applied Mathematics, Georg-August-University G\"ottingen, D-37083 G\"ottingen, Germany}

\begin{abstract}

This paper presents heavily grad-div and pressure jump stabilised, equal- and mixed-order discontinuous Galerkin finite element methods for non-isothermal incompressible flows based on the Oberbeck--Boussinesq approximation. In this framework, the enthalpy-porosity model for multiphase flow in melting and solidification problems can be employed. By considering the differentially heated cavity and the melting of pure gallium in a rectangular enclosure, it is shown that both boundary layers and sharp moving interior layers can be handled naturally by the proposed class of non-conforming methods. Due to the stabilising effect of the grad-div term and the robustness of discontinuous Galerkin methods, it is possible to solve the underlying problems accurately on coarse, non-adapted meshes. The interaction of heavy grad-div stabilisation and discontinuous Galerkin methods significantly improves the mass conservation properties and the overall accuracy of the numerical scheme which is observed for the first time. Hence, it is inferred that stabilised discontinuous Galerkin methods are highly robust as well as computationally efficient numerical methods to deal with natural convection problems arising in incompressible computational thermo-fluid dynamics.  

\end{abstract}


Publisher's version: DOI \url{https://doi.org/10.1016/j.jcp.2017.01.055}  \\
\copyright~2017. This manuscript version is made available under the CC BY-NC-ND 4.0 license: \\
\url{https://creativecommons.org/licenses/by-nc-nd/4.0/}

\end{frontmatter}

\section{Introduction}	\label{sec:Introduction}

In this paper, we propose a new class of stabilised dG-FEM for the numerical solution of the important category of incompressible computational thermo-fluid dynamic problems where the motion in the particular fluid is induced by natural convection phenomena due to local temperature differences. More precisely we are dealing with weakly non-isothermal flows following the Oberbeck--Boussinesq approximation \cite{Tritton88,KaysCrawford93}: 
\begin{subequations} \label{eq:OBEQs}
	\begin{empheq}[left=\empheqlbrace]{align} 
	\frac{\partial\uu}{\partial t}+\rb{\uu\cdot\nabla}\uu+\nabla p&=\nu\Delta \uu+\beta\rb{T-T_\refrm}\g \\
	\nabla\cdot\uu &=0 \\
	\frac{\partial T}{\partial t}+\uu\cdot\nabla T &=\alpha\Delta T
	\end{empheq} 
\end{subequations}
This well-known mathematical model describes buoyancy-driven flows which, for example, occur in indoor airflow simulations \cite{LubeEtAl08}. In this context, local mass conservation is very important for energy balance considerations. For a more general discussion of the Oberbeck-Boussinesq model we refer the reader to \cite{Zeytounian03}. \\

It is well-known that for this kind of fluid flows, the poor mass conservation properties of standard conforming FEM may lead to a loss of accuracy in the approximated solution \cite{GalvinEtAl12,JohnEtAl16,LinkeMerdon16}. Even though there are many conceivable remedies for this problem as, for example, the use of exactly divergence-free schemes \cite{CaseEtAl11,CockburnEtAl07,WangEtAl09}, or more generally the concept of pressure-robust methods \cite{Linke14,LinkeEtAl16}, we decided to take a different and original approach here. Namely, based on the idea of improving the fulfilment of the divergence constraint by means of grad-div stabilisation \cite{OlshanskiiEtAl09,JenkinsEtAl14}, we equip standard symmetric interior penalty dG-FEM \cite{PietroErn12,Riviere08,RoosEtAl08} with an additional grad-div term which can be implemented easily in any existing incompressible CFD code. Especially the natural treatment of convection-dominated problems, the inherent local mass conservation properties due to discontinuous pressures and the computational efficiency make dG-FEM appealing for the simulation of incompressible flows. In this context we consider both equal-order and mixed-order interpolation for velocity and pressure and therefore additionally introduce a pressure jump stabilisation term which is necessary for ensuring stability for the equal-order, but only optional for the mixed-order method \cite{SchoetzauEtAl03,CockburnEtAl09}.\\

This class of dG-FEM is applied to solve the following two different thermo-fluid flows:
\begin{itemize}
\item \textbf{Differentially heated cavity}. We consider Rayleigh numbers ranging from \num{e4} to \num{e8} and therefore have to deal with both velocity and thermal boundary layers of different thickness. Our solutions are analysed and compared to high-accuracy reference solutions, standard conforming FEM and an exactly divergence-free $\Hdiv$-conforming method. This analysis can be considered as the main part of this work.
\item \textbf{Melting of pure gallium}. Due to the occurring phase transition in this problem we have to deal with a solid/liquid phase interface which embodies a moving interior layer. A mesh sensitivity analysis and a parametric study for different sharpnesses of the interior layer is provided.
\end{itemize}

Both problems are solved on non-adapted, uniform fixed meshes which represents an additional numerical difficulty. In this work, however, we show that stabilised dG-FEM can be used successfully on non-adapted meshes\textemdash{}a suitable local mesh adaption would only improve our results. Due to the robustness of the proposed fully non-conforming schemes in combination with the enhanced discrete mass conservation properties, both boundary layer and moving interface problems can be handled easily. To the authors' knowledge, this work constitutes the first publication which combines dG-FEM with (heavy) grad-div stabilisation. \\

This paper is organised as follows. The mathematical modelling of one- and two-phase natural convection flows is given in Section \ref{sec:Model} where the enthalpy-porosity method is explained. The model for the heated cavity follows directly as a particular simplification. In order to obtain approximate solutions, the numerical treatment via stabilised dG-FEM is introduced in Section \ref{sec:Numerics}. Applying this class of dG-FEM, in Section \ref{sec:DHC}, we consider the differentially heated square cavity, and the melting of pure gallium is considered in Section \ref{sec:Gallium}. Finally, we summarise the results and conclude this work in Section \ref{sec:Conclusions}.

\section{One- and two-phase Oberbeck--Boussinesq model}	\label{sec:Model}
It turns out that the Oberbeck--Boussinesq model \eqref{eq:OBEQs} can be seen as a one-phase variant of the two-phase enthalpy-porosity model, applicable for fluid flows with solid/liquid phase transition phenomena \cite{BrentEtAl88,VollerEtAl89,Nikrityuk11,BelhamadiaEtAl12}. The development of numerical methods for such melting and solidification problems is numerically demanding due to the presence of a moving solid/liquid phase interface which embodies a sharp interior layer within the considered domain. \\

From a theoretical point of view, the main advantage of the enthalpy-porosity method is that the position of the moving solid/liquid interface is computed implicitly, solely depending on the local temperature. Thereby, conservation of mass, momentum and energy across the interface is ensured automatically. From a computational point of view the fixed-grid enthalpy-porosity method circumvents the need of developing more complicated numerical schemes which, for example, capture the phase boundary explicitly by means of adaptive mesh refinement. Comprehensive reviews, also of different approaches for modelling solid/liquid phase change problems, can be found, for example, in \cite{SamarskiiEtAl93,HuArgyropoulos96,Voller96}. In particular, we refer the reader to \cite{BelhamadiaEtAl04,DanailaEtAl14} for the use of adaptive meshes, to \cite{JanaEtAl07} for the employment of moving grids and to \cite{ChessaEtAl02} for an extended FEM in the context of the numerical solution of melting and solidification problems. For the treatment of the general case of free-surface and free-boundary flows, we recommend the reviews \cite{Caboussat05,ElgetiSauerland16}.\\

Typical for the enthalpy-porosity method is the introduction of a (smooth) phase change variable $\phi\colon \R\to\sqb{0,1}$, which depends on the local temperature and indicates solid and liquid phase:
\begin{equation}
\phi =
\begin{cases}
0,&\text{in the solid phase}\\
1,&\text{in the liquid phase}
\end{cases}
\end{equation}
Supposing that we can distinguish sharply between both phases, $\phi$ is discontinuous at the temperature of fusion which corresponds to a so-called `isothermal phase change' and physically means that during the melting process, the enthalpy increases by the latent heat without causing an increase in temperature. Typically, an isothermal phase change takes place for the melting and solidification of pure substances which have a precisely accurate and sharp melting temperature. For any other substances though, instead we have a certain non-zero temperature bandwidth, the so-called `mushy region', over which phase change occurs\textemdash{}these are so-called `non-isothermal phase change' processes. \\

\begin{figure}[h]
\centering
\includegraphics[width=0.6\textwidth]{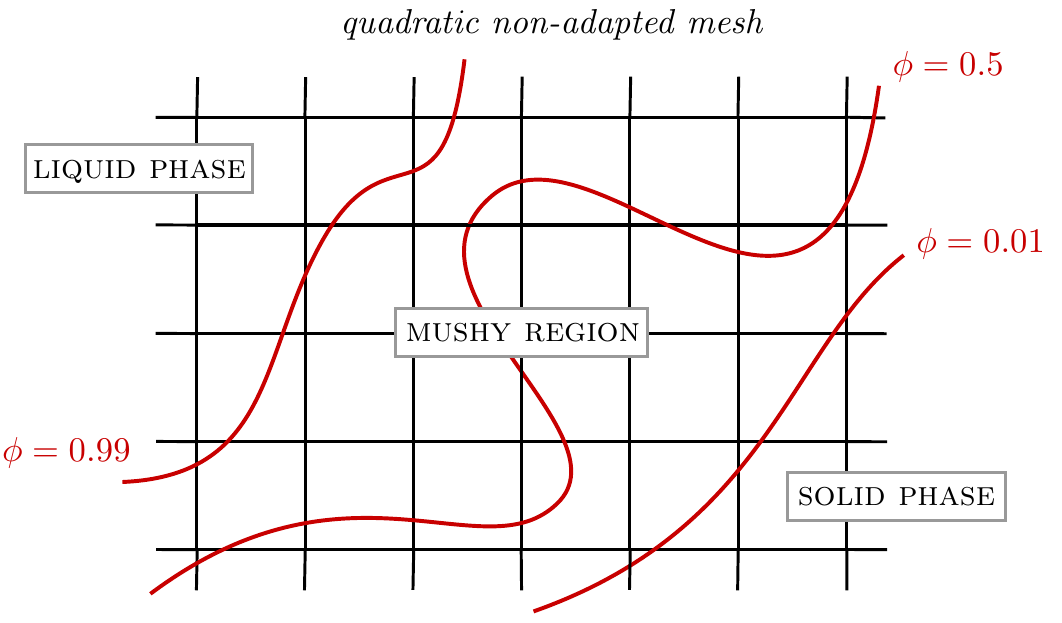}
\caption{Cutout of a quadratic non-adapted mesh (black lines) used for enthalpy-porosity methods. The $0.01$-, $0.5$- and $0.99$-contours (red lines) of the phase indicator function $\phi$ are shown which intersect the mesh elements arbitrarily. Dependent on $\phi$, the solid ($\phi\leqslant 0.01$) and liquid ($\phi\geqslant 0.99$) phase, and the mushy region ($0.01<\phi < 0.99$) is indicated. }
\label{fig:1}
\end{figure}

As pointed out in \cite{JanaEtAl07}, with regard to the validity of the enthalpy-porosity model, we want to emphasise that it is particularly powerful and straightforward for non-isothermal phase change processes. However, when there is need to obtain a phase change interface with zero thickness, the enthalpy-porosity method reaches its limit mainly because it is usually employed with a non-adapted mesh. To illustrate this problem, Figure \ref{fig:1} shows the moving phase boundary, represented by certain contours of the phase change indicator $\phi$, which typically intersect the mesh elements of a non-adapted mesh in an arbitrarily complex fashion. \\

Supposing that the liquid phase of the considered material behaves like a Newtonian fluid subject to incompressible, laminar flow with a constant kinematic viscosity $\nu$, the Navier--Stokes equations and the continuity equation are to be solved for obtaining the velocity $\uu$ and kinematic pressure $p$. Neglecting viscous dissipation, thermal radiation and both adiabatic compression and expansion, the standard energy equation with an additional source/sink term is used for obtaining the temperature $T$ in the material which is supposed to have a constant thermal conductivity $\kappa$. When the density $\rho$ and specific heat capacity at constant pressure $c_p$ are constant and equal in both the liquid and solid phase of the material, we thus obtain a constant thermal diffusivity given by
\begin{equation}
\alpha=\frac{\kappa}{\rho c_p}.	
\end{equation}
A straightforward extension of the enthalpy-porosity method that allows different thermodynamical properties between the two phases by using an effective volumetric heat capacity and thermal conductivity is proposed in \cite{BelhamadiaEtAl12}. Lastly, denote by $T_f$ the temperature of fusion of the material and let $\Delta T_f$ be the temperature range representing the width of the corresponding mushy region. In order to ensure differentiability, we define the above explained phase change indicator $\phi$ by the smooth hyperbolic tangent \cite{CagnoneEtAl14}
\begin{equation}
	\phi=F_{\Delta T_f}\rb{T}=
	\frac{1}{2}\sqb{\tanh\rb{\frac{5\rb{T-T_f}}{\Delta T_f}}+1},
\end{equation}
which is supposed to equal unity in the liquid phase and vanish in the solid phase. \\

Let $\Omega\subset \R^d$ be a bounded domain for $d\in\left\{2,3\right\}$ with Lipschitz boundary $\partial\Omega$ and $\tend >0$ the final time considered in the particular problem. The strongly coupled  nonlinear set of PDEs representing the enthalpy-porosity model for solid/liquid phase change problems reads as follows \cite{BrentEtAl88,VollerEtAl89,VollerEtAl87,VollerPrakash87}:
\begin{subequations} \label{eq:EnthalpyPorosityEQs}
	\begin{empheq}[left=\empheqlbrace]{alignat=2} 
	\frac{\partial\uu}{\partial t}+\rb{\uu\cdot\nabla}\uu+\nabla p&=\nu\Delta \uu+A\rb{\phi}\uu+\ff_b\rb{T},
	\qquad & \rb{\x,t}\in\Omega\times\left(0,\tend\right] \label{eq:EPNSEs}\\
	\nabla\cdot\uu&=0, & \rb{\x,t}\in\Omega\times\sqb{0,\tend} \label{eq:EPConti}\\
	\rho c_p\frac{\partial T}{\partial t}+\rho c_p\rb{\uu\cdot\nabla T} 
	&=\kappa\Delta T+\Q\rb{\phi}, & \rb{\x,t}\in\Omega\times\left(0,\tend\right] \label{eq:EPEnergy}\\
	\phi &= F_{\Delta T_f}\rb{T}, & \rb{\x,t}\in\Omega\times\sqb{0,\tend} 
	\end{empheq} 
\end{subequations}
We assume that the Oberbeck--Boussinesq approximation for weakly non-isothermal incompressible flows is valid \cite{GrayGiorgini76}, thereby yielding the general volumetric source term 
\begin{equation}
	\ff_b\rb{T}=\beta\rb{T-T_\refrm}\g
\end{equation}
in the Navier--Stokes equations which induces buoyancy effects due to the presence of gravitational forces. Here, $\beta$ denotes the coefficient of thermal expansion, $T_\refrm$ a reference temperature and $\g$ the vector representing the gravitational acceleration. The interface source term responsible for the absorption/release of energy during melting/solidification of a non-isothermal solid/liquid phase change process is given by \cite{BrentEtAl88,BelhamadiaEtAl12}
\begin{equation}
	\Q\rb{\phi}=-\rho L_f\sqb{\frac{\partial \phi}{\partial t}+\uu\cdot\nabla \phi},
\end{equation}
where we note that this term only acts in the  mushy region. Here, $L_f$ denotes the latent heat of fusion. Lastly, the term which achieves that the solid material has zero velocity is chosen as
\begin{equation}
	A\rb{\phi}=-\frac{C_0}{\rho}\frac{\rb{1-\phi}^2}{|\phi|^3+b},
\end{equation}
where $C_0> 0$ denotes a large parameter responsible for the attenuation both in the mushy region and the solid phase and $b> 0$ is a security parameter preventing division by zero whenever $\phi\equiv 0$. This term is inspired by the Carman--Kozeny equations \cite{VollerEtAl89,VollerPrakash87,Carman97}. We notice that the original one-phase Oberbeck--Boussinesq model \eqref{eq:OBEQs} is a special case of the two-phase enthalpy-porosity model \eqref{eq:EnthalpyPorosityEQs} with $\phi\equiv 1$. \\ 

To close system \eqref{eq:EnthalpyPorosityEQs}, initial and boundary values are required. Therefore, we assume that the boundary $\partial\Omega$ can be decomposed into two pairwise disjoint sets $\Gamma_D$ and $\Gamma_N$ for which $\partial\Omega=\Gamma_D\cup\Gamma_N$ holds true. On the Dirichlet part $\Gamma_D$ we prescribe the temperature $T$ as a known, possibly time-dependent, function $g_D^T$. On the other hand, a prescribed heat flux $g_N^T$ across the boundary is specified for the Neumann part $\Gamma_N$, i.e.
\begin{equation}
	\kappa\nabla T\cdot \n=g_N^T\quad\text{on}~\Gamma_N.
\end{equation}
As usual, $\n$ denotes the outward unit normal vector to $\partial\Omega$. With respect to the velocity field the no-slip condition $\uu=\mathbf{0}$ is assumed to hold true on $\partial\Omega$, and for the pressure we impose the zero-mean condition. Initially, the temperature in the whole domain $\Omega$ is known and defined by the function $T_0$ for which compatibility with the boundary conditions has to be ensured. In all subsequent problems, the fluid is initially at rest, hence $\uu\equiv \mathbf{0}$ and  $p\equiv 0$ at $t=0$. 

\section{Stabilised dG-FEM} \label{sec:Numerics}

In order to find approximate solutions to \eqref{eq:EnthalpyPorosityEQs} this section proposes differently stabilised dG-FEM based on both mixed- and equal-order interpolation for the velocity and pressure. Therefore, basic notations for the treatment of dG-FEM are introduced which enable the statement of the variational formulation of the generic problem. In this formulation, several numerical parameters occur which are clarified afterwards. For more details concerning standard dG-FEM we refer the reader to \cite{PietroErn12,Riviere08,RoosEtAl08}, on which the numerical part of this work is loosely based. The following explanations are, for the sake of simplicity, restricted to the spatially two-dimensional case. Nonetheless, an extension to the three-dimensional case is straightforward. 

\subsection{Preliminaries}

Suppose that the polygonal domain $\Omega\subset \R^2$ is partitioned into an admissible and quasi-uniform decomposition consisting of quadrilateral mesh elements $\T=\left\{K_1,\hdots,K_M\right\}$ such that $\overline{\Omega}=\cup_{K\in\T}\overline{K}$. The subscript $h$ refers to the refinement of the mesh and is defined by 
\begin{equation}
	h=\max_{K\in\T} h_K,\quad h_K=\max_{F\subset\partial K} h_F,
\end{equation}
where $h_F$ denotes the diameter of the edge $F$. For any element $K\in\T$ its outward unit normal vector is given by $\n_K$. Moreover, let $\F$ denote the set of all edges corresponding to $\T$ and $\FK=\left\{F\in\F\colon F\subset\partial K\right\}$. Let $\Fi\subset\F$ be the subset of all interior edges and $\Fb\subset\F$ denote the subset of all boundary edges.  To any edge $F\in\F$ we assign a unit normal vector $\n_F$ where for edges $F\in\Fb$ this is the usual outward unit normal vector $\n$. If $F\in\Fi$, there are two elements $K_1^F$ and $K_2^F$ sharing the edge $F$ and $\n_F$ is supposed to point in an arbitrary, but fixed, direction. On the interface between two neighbouring elements, the average and jump operators are defined for any piecewise continuous function $v$:
\begin{subequations}
\begin{empheq}{alignat=3} 
&\text{Average:}\quad &&\avg{v}_F=\frac{1}{2}\rb{\restr{v}{{K_1^F}}+\restr{v}{{K_2^F}}},\quad &&\forall F=\partial{K_1^F}\cap\partial{K_2^F}\in \Fi \\
&\text{Jump:} &&\jmp{v}_F= \restr{v}{{K_1^F}}-\restr{v}{{K_2^F}},\quad &&\forall F=\partial{K_1^F}\cap\partial{K_2^F}\in\Fi
\end{empheq}	
\end{subequations}  
For boundary edges $F\in\Fb$ we set $\avg{v}_F=\jmp{v}_F=\restr{v}{F}$. Furthermore, when dealing with vector-valued functions both the average and jump operator are supposed to act componentwise. When no confusion can arise the subscript indicating the edge is omitted and we simply write $\avg{\cdot}$ and $\jmp{\cdot}$. \\

For $k\geqslant 1$ we define the discontinuous finite element space
\begin{equation}
	\V_h^k=\left\{v\in\Lp^2\rb{\Omega}\colon\restr{v}{K}\in\mathcal{Q}_k\rb{K},~\forall K\in\T\right\},
\end{equation}
where $\mathcal{Q}_k$ denotes the space of tensor product polynomials of degree $k$. We also refer to this  space as $\Quad{-k}$ and denote by $\Quad{k}$ its continuous counterpart. Additionally, define the following discontinuous FE spaces:
\begin{equation}
\UU_h^k=\sqb{\V_h^k}^2,\quad \PP_h^\kpr=\left\{q\in\V_h^\kpr\colon \int_\Omega q\rb{\x}\,\drm\x=0\right\}
\end{equation}
Note that for $\kpr=k-1$ we obtain a mixed-order Taylor--Hood type method whereas $\kpr=k$ yields a method with equal-order interpolation for velocity and pressure. The spaces $\V_h^k$, $\UU_h^k$ and $\PP_h^\kpr$ are used to state the variational formulation arising from a dG-FEM in search for approximate solutions $T_h$, $\uu_h=\rb{u_{1h},u_{2h}}^\dag$ and $p_h$ for the temperature $T$, velocity field $\uu$ and pressure $p$, respectively.

\subsection{Variational formulation}

In this work we concentrate on the spatial semi-discretisation of the enthalpy-porosity model \eqref{eq:EnthalpyPorosityEQs} by using differently stabilised dG-FEM. As usual, $\rb{\cdot,\cdot}$ denotes the $\Lp^2$-inner product on the whole domain $\Omega$. The semi-discrete variational formulation of this problem reads as follows:
\begin{subequations} \label{eq:DGVarForm}
\begin{empheq}[left=\empheqlbrace]{alignat=2} 
\text{find }\rb{\uu_h,p_h,T_h}\in\Lp^2\rb{\sqb{0,\tend};\UU_h^k\times \PP_h^\kpr\times \V_h^k}~\text{s.t.}~\forall\rb{\vv_h,q_h,v_h}&\in\UU_h^k\times \PP_h^\kpr\times \V_h^k\\
\rb{\partial_t \uu_h,\vv_h}+a_h\rb{\nu;\uu_h,\vv_h}+t_h\rb{\uu_h,\uu_h,\vv_h}+b_h\rb{\vv_h,p_h}	+j_h\rb{\uu_h,\vv_h}
&=\ell_h\rb{\uu_h,\vv_h} \label{eq:DGNSEs}\\
-b_h\rb{\uu_h,q_h}+s_h\rb{p_h,q_h}									
&=0 \label{eq:DGConti}\\
\rb{\rho c_p\partial_t T_h,v_h}+a_h^\sip\rb{\kappa;T_h,v_h}+a_h^\upw\rb{\rho c_p T_h,v_h}&
=\ell_h^\sip\rb{\kappa;v_h}+\ell_h^\upw\rb{v_h} \label{eq:DGEnergy}\\
\phi_h &= F_{\Delta T_f}\rb{T_h}
\end{empheq}
\end{subequations}
Here, \eqref{eq:DGNSEs} represents the Navier--Stokes equations solving for the discrete velocity $\uu_h$, \eqref{eq:DGConti} the continuity equation solving for the discrete kinematic pressure $p_h$ and \eqref{eq:DGEnergy} the energy equation solving for the discrete temperature $T_h$. We specify the linear, bilinear and trilinear forms occurring in problem \eqref{eq:DGVarForm} below. All numerical schemes in this work were implemented into the finite element package COMSOL Multiphysics 5.1 by the first author. For the time stepping in \eqref{eq:DGVarForm}, the fully implicit, second-order, variable-time-step BDF(2) method IDA from the SUNDIALS suite is employed \cite{HindmarshEtAl05,BrownEtAl94}. All occurring nonlinearities are treated by a Newton method with tolerance \num{e-10} where the corresponding linear systems are solved directly by the PARDISO solver which uses efficient parallel sparse LU factorisation \cite{SchenkEtAl00,Schenk00,SchenkGartner04}. The number of Newton iterations per time step is restricted to 6 and the Jacobian is updated only once per time step. Since our BDF(2) method is adaptive in time, only a maximum time step can be imposed, which is chosen as $\Delta t_\maxrm=1$ for the differentially heated cavity and $\Delta t_\maxrm=0.05$ for the melting of pure gallium. Therefore, the particular time step size depends on the fulfilment of the tolerance criterion. In order to give an impression for this, on average, for the heated cavity our simulations need $\approx{1000}$ total time steps to finish whilst for the gallium problem $\approx{2500}$ time steps are necessary.

\subsubsection{Discrete Laplacian and source terms}

In order to account for the Laplace operators in \eqref{eq:EPNSEs} and \eqref{eq:EPEnergy} we use a symmetric interior penalty (SIP) formulation of the dG-FEM \cite{PietroErn12,BottiPietro11}. The consistent, symmetric and parameter-dependent SIP bilinear form $a_h^\sip\colon \V_h^k\times \V_h^k\to\R$ is defined as
\begin{equation}
	\begin{aligned}
a_h^\sip\rb{\eps;w_h,v_h}
	=&	\int_\Omega \eps\nabla_h w_h \cdot \nabla_h v_h\,\drm\x
		-\sum_{F\in\FbD}\oint_F \sqb{\rb{\eps\nabla w_h\cdot \n}v_h 
		+ w_h \rb{\eps\nabla v_h\cdot \n}
		-\eta\frac{\eps}{h_F} w_h v_h}\,\drm\s\\
	-&	\sum_{F\in\Fi}\oint_F \sqb{\avg{\eps\nabla w_h}\cdot \n_F \jmp{v_h}
		+ \jmp{w_h}\avg{\eps\nabla v_h}\cdot \n_F
		-\eta\frac{\eps}{h_F}\jmp{w_h}\jmp{v_h}}\,\drm\s	,
\end{aligned} 
\end{equation} 
where the generic diffusion parameter $\eps>0$ stands for either the thermal conductivity $\kappa$ or the kinematic viscosity $\nu$. Note that the operator $\nabla_h$ denotes the broken gradient which is defined for any piecewise continuous function $v$ by
\begin{equation}
 	\restr{\rb{\nabla_h v}}{K}=\nabla\rb{\restr{v}{K}},\quad \forall K\in\T.
 \end{equation}
It is well-known that the SIP bilinear form is bounded and enjoys discrete coercivity whenever the discontinuity penalisation parameter $\eta>0$ is sufficiently large \cite{PietroErn12,Riviere08}. The particular choice of this parameter is discussed in Section \ref{sec:DiscPena}. \\

Both homogeneous and non-homogeneous Dirichlet boundary conditions are enforced weakly following Nitsche's method \cite{Arnold82}. Therefore, for approximating the temperature in \eqref{eq:DGEnergy}, the likewise parameter-dependent linear form $\ell_h^\sip\colon \V_h^k\to \R$ accounts for the corresponding source term and non-homogeneous boundary conditions for the temperature \cite{PietroErn12,Riviere08}
\begin{align}
	\ell_h^\sip\rb{\eps;v_h}
	=&	\int_\Omega \Q\rb{\phi} v_h \,\drm\x 
		-\sum_{F\in\FbD}\oint_F \sqb{\rb{\eps\nabla v_h\cdot \n}- \eta\frac{\eps}{h_F} v_h}g_D^T \,\drm\s
		+ \sum_{F\in\FbN}\oint_F g_N^T v_h \,\drm\s.
\end{align} 
Here, $g_D^T$ denotes the Dirichlet boundary condition and $g_N^T$ represents the Neumann data on the boundary edges $\FbD=\Fb\cap\Gamma_D$, belonging to the Dirichlet boundary, and $\FbN=\Fb\cap\Gamma_N$, belonging to the Neumann boundary, respectively. 	The occurring time derivative of the phase change indicator $\phi$ in the interface source term $\Q\rb{\phi}$ is treated in a fully implicit manner. That is, we obtain $\frac{\partial \phi}{\partial t}=\frac{\partial \phi}{\partial T}\frac{\partial T}{\partial t}$ by the chain rule since $\phi=\phi\rb{T}$ is implicitly defined over the local temperature. \\

For approximating the velocity field in \eqref{eq:DGNSEs} the SIP bilinear form is employed, as well. Supposing that $\ww_h=\rb{w_{1h},w_{2h}}^\dag$ and $\vv_h=\rb{v_{1h},v_{2h}}^\dag$, the parameter-dependent bilinear form $a_h\colon\UU_h^k\times\UU_h^k\to\R$ applies the SIP bilinear form componentwise to the components of the velocity as follows \cite{BottiPietro11,PietroErn10}:
\begin{align}
	a_h\rb{\nu;\ww_h,\vv_h}=a_h^\sip\rb{\nu;w_{1h},v_{1h}}+a_h^\sip\rb{\nu;w_{2h},v_{2h}}
\end{align}
Since the no-slip condition is imposed, there is no need to account for non-homogeneous boundary conditions, thereby simplifying the corresponding form $\ell_h\colon\UU_h^k\times\UU_h^k\to\R$ for the velocity to 
\begin{align}
\ell_h\rb{\uu_h,\vv_h}=\int_\Omega \sqb{A\rb{\phi_h}\uu_h\cdot\vv_h + \ff_b\rb{T_h}\cdot\vv_h}\,\drm\x.
\end{align}
Note that both the velocity attenuation reaction term and the Oberbeck--Boussinesq source term are condensed in this form.
\subsubsection{Discrete convective forms}

Inherent to computational thermo-fluid dynamics is the presence of convective heat and mass transfer. Whereas in the Navier--Stokes equations \eqref{eq:EPNSEs} this is represented by the nonlinear inertia term $\rb{\uu\cdot\nabla}\uu$, the energy equation \eqref{eq:EPEnergy} comprises the convective term $\uu\cdot\nabla T$. In order to account for the latter, the following upwind bilinear form $a_h^\upw\colon \V_h^k\times \V_h^k\to\R$ is employed for computing the approximated temperature \cite{PietroErn12,ErnEtAl09,BrezziEtAl04}:
\begin{equation} \label{eq:TempUpw}
\begin{aligned}
a_h^\upw\rb{w_h,v_h}
	=&	\int_\Omega \rb{\uu_h\cdot\nabla_h w_h} v_h\,\drm\x
		+ \sum_{F\in\FbD}\oint_F \rb{\uu_h\cdot\n}^\ominus w_hv_h\,\drm\s\\
	-&	\sum_{F\in\Fi}\oint_F\sqb{\rb{\avg{\uu_h}\cdot\n_F}\jmp{w_h}\avg{v_h} 
		- \frac{1}{2}\abs{\avg{\uu_h}\cdot \n_F}\jmp{w_h}\jmp{v_h}}\,\drm\s
\end{aligned}	
\end{equation}
In doing so, the negative part $v^\ominus=\frac{1}{2}\rb{\abs{v}-v}$ of a function $v$ is used. For non-homogeneous Dirichlet boundary conditions for the temperature, the upwind bilinear $a_h^\upw\rb{\cdot,\cdot}$ form cooperates with the linear form $\ell_h^\upw\colon \V_h^k\to \R$ defined by 
\begin{align}
	\ell_h^\upw\rb{v_h}=\sum_{F\in\FbD}\oint_F \rb{\uu_h\cdot\n}^\ominus g_D^T v_h\,\drm\s.
\end{align}
The inertia term in the Navier--Stokes equations is treated by the following trilinear form $t_h\colon\UU_h^k\times\UU_h^k\times\UU_h^k\to\R$ which incorporates Temam's modification on the discrete level \cite{BottiPietro11,PietroErn10}:
\begin{equation}\label{eq:InertiaTerm}
\begin{aligned}
	 t_h\rb{\ww_h,\uu_h,\vv_h}
	=&	\int_\Omega \sqb{\rb{\ww_h\cdot\nabla_h}\uu_h\cdot\vv_h
		+\frac{1}{2}\rb{\nabla_h\cdot\ww_h}\rb{\uu_h\cdot\vv_h}}\,\drm\x\\
	-&	\sum_{F\in\Fi}\oint_F \sqb{\rb{\avg{\ww_h}\cdot\n_F}\rb{\jmp{\uu_h}\cdot\avg{\vv_h}}
		+\frac{1}{2}\rb{\jmp{\ww_h}\cdot\n_F}\avg{\uu_h\cdot\vv_h}}\,\drm\s\\
	-&	\frac{1}{2}\sum_{F\in\Fb}\oint_F \rb{\ww_h\cdot\n}\rb{\uu_h\cdot\vv_h}\,\drm\s
\end{aligned}	
\end{equation}
Here, $\nabla_h\cdot$ denotes the broken divergence which acts elementwise. This trilinear form is not locally mass conservative but contains a term proportional to the broken divergence of the discrete velocity. In order to improve the mass conservation properties of the proposed method the following subsection introduces  grad-div stabilisation. \\

We briefly want to comment on our choice of convective fluxes. Surprisingly, based on our experience with the considered weakly non-isothermal flows in the laminar regime, it is not necessary to include any form of upwind stabilisation for the velocity approximation in \eqref{eq:InertiaTerm}. For the temperature approximation, however, we included a standard upwinding term in \eqref{eq:TempUpw} which, especially on coarse meshes, clearly improves the discrete solution.

\subsubsection{Discrete velocity-pressure coupling and stabilisation}
The pressure-velocity coupling, also called the discrete divergence, is realised with the bilinear form $b_h\colon\UU_h^k\times \PP_h^\kpr\to\R$ defined by \cite{BottiPietro11,PietroErn10,GiraultEtAl05}
\begin{align}
b_h\rb{\ww_h,q_h}
	=-\int_\Omega q_h\rb{\nabla_h\cdot\ww_h}\,\drm\x+\sum_{F\in\Fb}\oint_F \rb{\ww_h\cdot\n}q_h\,\drm\s
	+\sum_{F\in\Fi}\oint_F \rb{\jmp{\ww_h}\cdot\n} \avg{q_h} \,\drm\s.
\end{align}
In order to ensure the stability of the method with equal-order $\rb{\Quad{-2}\slash\Quad{-2}}$ interpolation for velocity and pressure the bilinear form $s_h\colon\PP_h^\kpr\times\PP_h^\kpr\to\R$, which penalises pressure jumps across interfaces, is introduced \cite{PietroErn10,CockburnEtAl09}:
\begin{align}
s_h\rb{q_h,r_h}
	= \sum_{F\in\Fi}\lambda \frac{h_F}{\nu}\oint_F\jmp{q_h}\jmp{r_h}\,\drm\s
\end{align}
For Taylor--Hood type $\rb{\Quad{-2}\slash\Quad{-1}}$ dG-FEM pressure jump stabilisation is not necessary to ensure stability \cite{Riviere08}. Nevertheless, originating from local discontinuous Galerkin (LDG) methods \cite{CockburnEtAl05a}, such a term can also be included successfully in the mixed-order formulation \cite{SchoetzauEtAl03}. Therefore, basically three different methods can be obtained which we want to compare in the following:
\begin{enumerate}[label=(\roman*)]
\item Taylor--Hood type dG-FEM ($\kpr=k-1$) without pressure jump penalisation ($\lambda=0$)
\item Taylor--Hood type dG-FEM ($\kpr=k-1$) with pressure jump penalisation ($\lambda>0$)
\item Equal-order type dG-FEM ($\kpr=k$) with pressure jump penalisation ($\lambda>0$)	
\end{enumerate}
Lastly, the mass conservation properties of all considered methods are to be improved. It is well-known that for coupled flow problems poor mass conservation, especially in conjunction with large and complex pressures, results in a loss of accuracy of the approximated solution \cite{GalvinEtAl12,Linke09}. One simple to implement possible remedy constitutes the grad-div stabilisation which stems from adding the consistent term $-\gamma\nabla\rb{\nabla\cdot\uu}$ with $\gamma\geqslant 0$ to the left-hand side of the Navier--Stokes equations \cite{OlshanskiiEtAl09,CaseEtAl11}. In the corresponding variational formulation the bilinear form $j_h\colon \UU_h^k\times\UU_h^k\to\R$ accounts for grad-div stabilisation:
\begin{align}
j_h\rb{\ww_h,\vv_h}
	=\sum_{K\in\T}\gamma\int_K\rb{\nabla\cdot\ww_h}\rb{\nabla\cdot\vv_h}\,\drm\x
\end{align}
It is our aim to show that it is possible to devise dG-FEM where heavy grad-div stabilisation, $\gamma\gg 1$, does not reduce the accuracy of the approximated solution. Furthermore, we will see that in such a way essentially pointwise divergence-free velocity fields can be obtained easily and thus, the mass conservation of the proposed method is significantly improved. Note that the idea of using heavy grad-div stabilisation is not novel; cf.\ \cite{GalvinEtAl12,JenkinsEtAl14} in the context of conforming FEM. For dG-FEM, however, to our knowledge this work is the first reported result which shows that a benefit can be obtained also for non-conforming methods.

\subsection{Discontinuity penalisation parameter} \label{sec:DiscPena}
In stating the variational formulation \eqref{eq:DGVarForm} there is one parameter which is not yet defined properly. In order to account for second-order spatial derivatives, corresponding to diffusion phenomena, the discontinuity penalisation parameter $\eta$ has been introduced. We note that due to this term, interior penalty methods aim at approximate continuity of the considered variables. For the vector-valued velocity this means that the penalty suppresses singularities of discrete gradient and divergence on the inter-element boundaries. \\

If $\eta$ is chosen too small, discrete coercivity of the formulation cannot be guaranteed \cite{PietroErn12}. As pointed out in \cite{PietroErn12,Riviere08} there exists a minimum penalty parameter $\eta^*$, dependent on the maximum number $N_\partial$ of neighbours an element $K$ of the decomposition $\T$ can have and the constant $C_\tr$ in the discrete trace inequality
\begin{equation}
\forall F\in\FK:\quad
\norm{v_h}_{\Lp^2\rb{F}}\leqslant C_\tr h_K^{-\nicefrac{1}{2}}\norm{v_h}_{\Lp^2\rb{K}},
\quad \forall v_h\in \mathcal{Q}_k\rb{K}.
\end{equation}
This minimum penalty parameter is then given by \cite{PietroErn12}
\begin{equation}
\eta^*= C_\tr^2 N_\partial,	
\end{equation}
where for interior elements $N_\partial=4$ since, in this work, $\T$ consists of quadrilateral elements and hanging nodes are prohibited. For elements in contact with the boundary, we thus obtain $N_\partial=3$. In \cite{Hillewaert13} a sharp bound for $C_\tr$ is developed and proven for quadrilateral mesh elements in two space dimensions:
\begin{equation}
\forall F\in\FK:\quad\norm{v_h}_{\Lp^2\rb{F}}^2\leqslant \rb{k+1}^2\frac{\Area{F}}{\Vol{K}}\norm{v_h}_{\Lp^2\rb{K}}^2, \quad \forall v_h\in \mathcal{Q}_k\rb{K}
\end{equation}
Here, $\Area{F}$ and $\Vol{K}$ denote length and area for $F$ and $K$, respectively. In the following, only quadratic decompositions are considered for which $\frac{\Area{F}}{\Vol{K}}= h_K^{-1}$ holds true for all $F\in\FK$ and therefore
\begin{equation}
C_\tr^2=\rb{k+1}^2.	
\end{equation}
Furthermore, for PDEs involving diffusion phenomena, we restrict ourselves to exclusively using biquadratic, $k=2$, interpolation yielding $C_\tr^2=9$. The resulting minimum penalisation parameter is thus given by $\eta^*=36$ for interior edges and $\eta^*=27$ for boundary edges. For all subsequent simulations in this work, as the jump penalisation parameter we use $\eta=\eta^*$ for the velocity and $\eta=72>\eta^*$ for the temperature.

\section{Boundary layers: One-phase flow in a differentially heated cavity} \label{sec:DHC}

The differentially heated cavity (DHC) with adiabatic top and bottom walls is a classical problem considered in the literature on heat transfer processes with fluid flow \cite{Davis83,Quere91,BarakosEtAl94,XinQuere06}. Therefore, there is a wide variety of solutions available which makes this problem a very well-suited test case for assessing the quality of the proposed stabilised dG-FEM. The model of the DHC problem is given by the coupled system \eqref{eq:OBEQs} which is equivalent to the enthalpy-porosity system \eqref{eq:EnthalpyPorosityEQs} without phase change, i.e.\ $\phi\equiv 1$. Typically, we want to consider a dimensionless formulation of the problem and use the dimensionless numbers 
\begin{equation}
\Pra=\frac{\mu c_p}{\kappa}=\frac{\nu}{\alpha},\quad
\Gra= \frac{g\beta\Delta T_\refrm L_\refrm^3}{\nu^2}\quad\text{and}\quad 
\Ray= \Gra\,\Pra,
\end{equation}
where we restrict ourselves to a square cavity with height and width $L_\refrm=1$ and $\Delta T_\refrm>0$ denotes the difference in temperature of the two vertical walls. In this work as well as in all the above given references, the Prandtl number is fixed at
\begin{equation}
\Pra=0.71,
\end{equation}
corresponding to a cavity filled with air. Note that for Prandtl numbers below unity, velocity boundary layers are generally sharper than the corresponding thermal boundary layers \cite{Tritton88,Durst08}. For ease of notation, we do not introduce new symbols for the dimensionless quantities corresponding to the velocity, pressure, temperature, time and space domain. However, we keep in mind that all expressions in this section concerning the DHC problem do not have a physical unit. Following the formulation in \cite{Quere91}, the dimensionless governing equations in primitive, dimensionless variables, supposing that the dimensionless domain is given by the unit square $\Omega=\rb{0,1}^2$, are provided by \eqref{eq:EnthalpyPorosityEQs} with 
\begin{equation}
\nu=\frac{\Pra}{\sqrt{\Ray}},\quad
\phi \equiv \rho c_p\equiv 1,\quad
A\rb{\phi}\equiv\Q\rb{\phi} \equiv 0,\quad 
\ff_b\rb{T}=\rb{0,\Pra\,T}^\dag\quad \text{and}\quad 
\kappa=\alpha=\frac{1}{\sqrt{\Ray}}	.
\end{equation}
Note that the problem is formulated as a time-dependent system of PDEs but in fact, we are only interested in the steady-state solution, which is the solution for the time-dependent problem as $\tend\to\infty$. Note that such a steady-state solution, in two space dimensions, is known to become unstable for $\Ray\geqslant \num{2E8}$ \cite{QuereBehnia98}. \\

The following closing initial conditions (ICs) and boundary conditions (BCs) are imposed:
\begin{enumerate}[label=(\roman*)]
\item Initial values $u_1\equiv u_2 \equiv p \equiv 0$ and $T\equiv -\frac{1}{2}$ on $\Omega$ at $t=0$ and the no-slip condition $\uu \equiv \mathbf{0}$ on $\partial\Omega$.
\item Dirichlet conditions $g_D^T=T=\frac{1}{2}$ on $x_1=0$ and $g_D^T=T=-\frac{1}{2}$ on $x_1=1$ for all $0\leqslant x_2\leqslant 1$ for the temperature. Define by $\Gamma_D=\left\{\rb{x_1,x_2}\in\Omega\colon x_1=0~\text{and}~x_1=1\right\}$ the Dirichlet part of the boundary $\partial\Omega$ with prescribed Dirichlet boundary condition $g_D^T$ on $\Gamma_D$.
\item Homogeneous Neumann conditions $g_N^T=\frac{\partial T}{\partial x_2}=0$ on $x_2=0$ and $x_2=1$ for all $0\leqslant x_1\leqslant 1$ for the temperature. Define by $\Gamma_N=\left\{\rb{x_1,x_2}\in\Omega\colon x_2=0~\text{and}~x_2=1\right\}$ the Neumann part of the boundary $\partial\Omega$ with prescribed Neumann boundary condition $g_N^T$ on $\Gamma_N$.
\end{enumerate}

In order to ensure compatibility between initial and boundary conditions for the temperature, the temperature BC at $x_1=0$ is ramped up smoothly during the first ten time units of computation. Our simulations are computed for $\Ray\in\left\{\num{e4},\num{e6},\num{e8}\right\}$ on quadratic non-adapted meshes without any refinement towards the boundary. We will see that, even though especially for high Rayleigh numbers strong thermal and velocity boundary layers are forming, all proposed stabilised dG-FEM easily cope with under-resolved meshes. \\

\begin{figure}[h]
\centering
\subfigure[$\Ray=\num{e4}$; $\tend=100$]{\includegraphics[width=0.285\textwidth]{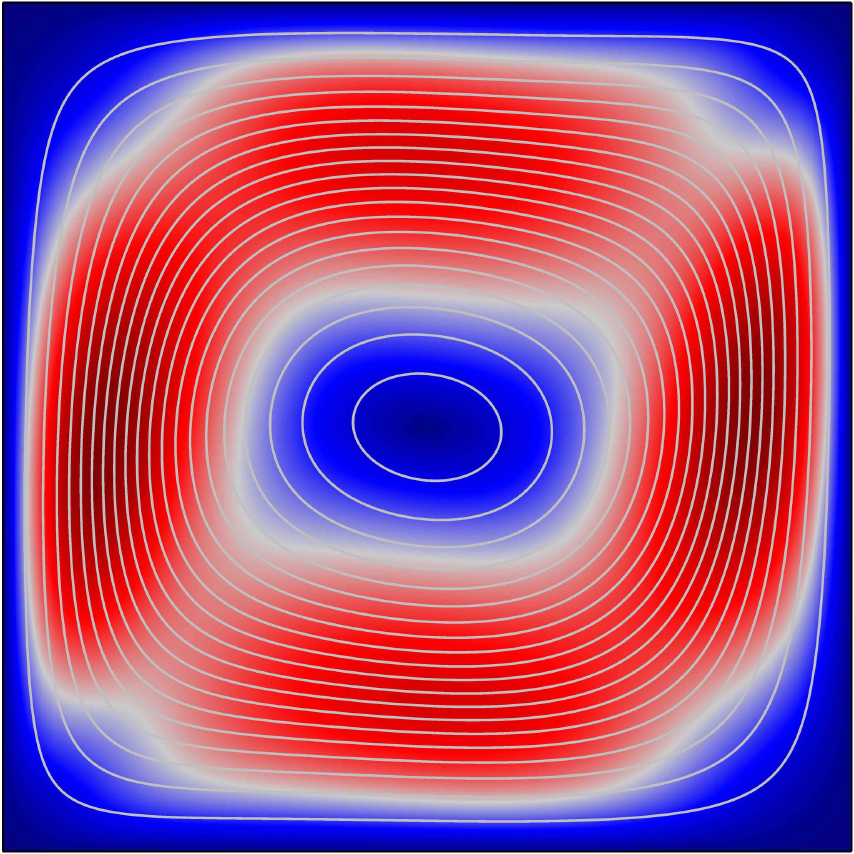}}\goodgap
\subfigure[$\Ray=\num{e6}$; $\tend=300$]{\includegraphics[width=0.285\textwidth]{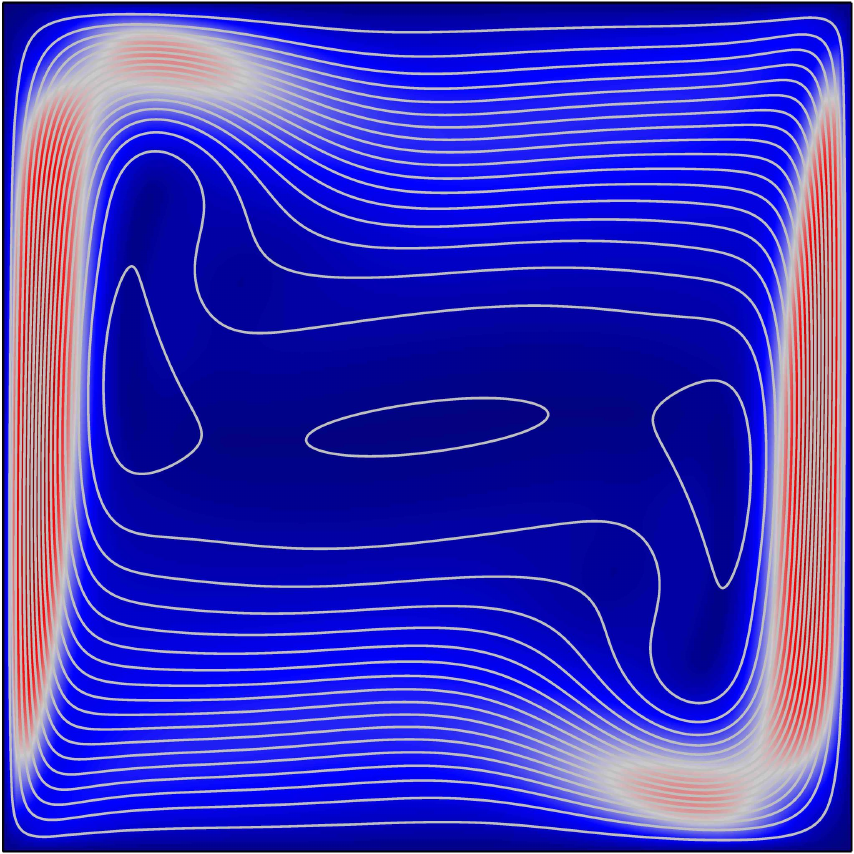}}\goodgap
\subfigure[$\Ray=\num{e8}$; $\tend=750$]{\includegraphics[width=0.285\textwidth]{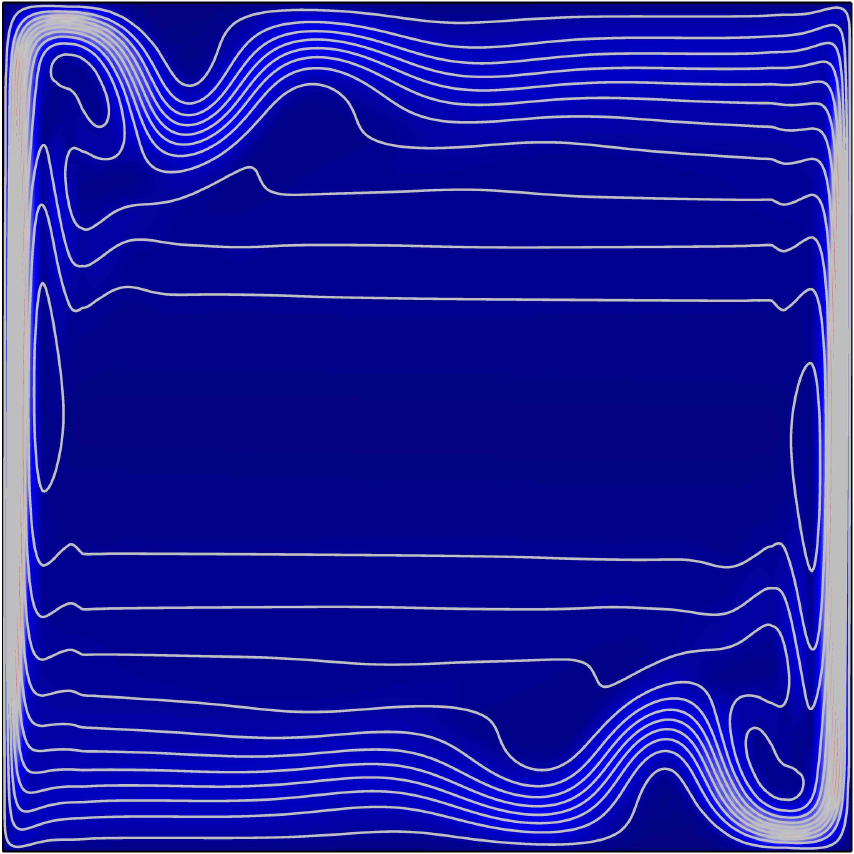}}
\caption{Velocity magnitude and streamlines of DHC simulation for $\Ray\in\left\{\num{e4},\num{e6},\num{e8}\right\}$. Results are computed with the $\rb{\Quad{-2}\slash\Quad{-1}}\wedge\Quad{-2}$ dG-FEM on $64\times 64$ meshes with stabilisation parameters $\gamma=\num{e5}$ and $\lambda=\num{e3}$.}
\label{fig:2}
\end{figure}

Figure \ref{fig:2} shows the streamlines of the velocity field obtained by the grad-div and pressure stabilised Taylor--Hood type $\rb{\Quad{-2}\slash\Quad{-1}}\wedge\Quad{-2}$ dG-FEM on $64\times 64$ meshes for different values of $\Ray$. First of all, we note that although we use a coarser mesh, the results are in excellent agreement with the above mentioned literature for the whole range of Rayleigh numbers considered. Analysing the flow, one can see that at the hot wall the fluid is rising whilst at the cold wall the fluid is dropping down, where converging streamlines indicate higher velocities. The behaviour of the flow in the middle of the cavity is strongly dependent on the particular Rayleigh number. For $\Ray=\num{e4}$ the flow field shows a central vortex with a slight tendency of becoming elliptic. As $\Ray$ increases to $\Ray=\num{e6}$ the vortex disassembles into two vortices which move to the vertical walls, thereby making room for a third small vortex in the centre. For $\Ray=\num{e8}$ the central vortex vanishes again and the two vortices move to the upper left and lower right corner, respectively. Apart from the vertical walls a stratified flow field is forming. Moreover, we note that the minimum end of time $\tend$ for the simulation which yields a stable stationary solution increases with $\Ray$. \\

\begin{figure}[h]
\centering
\subfigure[$\Ray=\num{e4}$; $\tend=100$]{\includegraphics[width=0.285\textwidth]{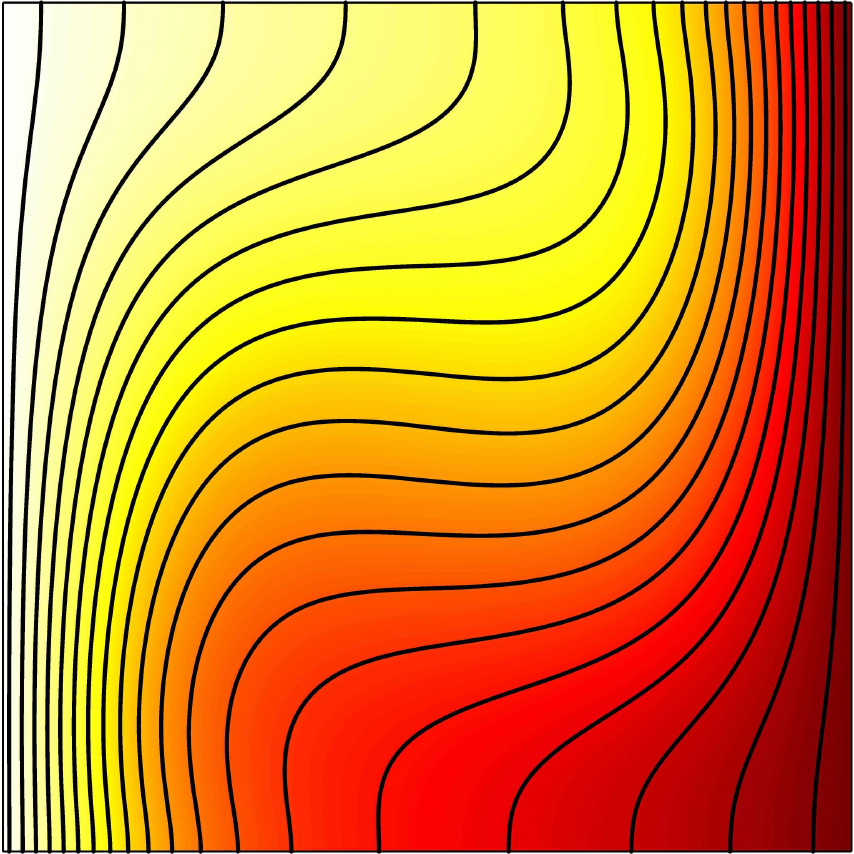}} \goodgap
\subfigure[$\Ray=\num{e6}$; $\tend=300$]{\includegraphics[width=0.285\textwidth]{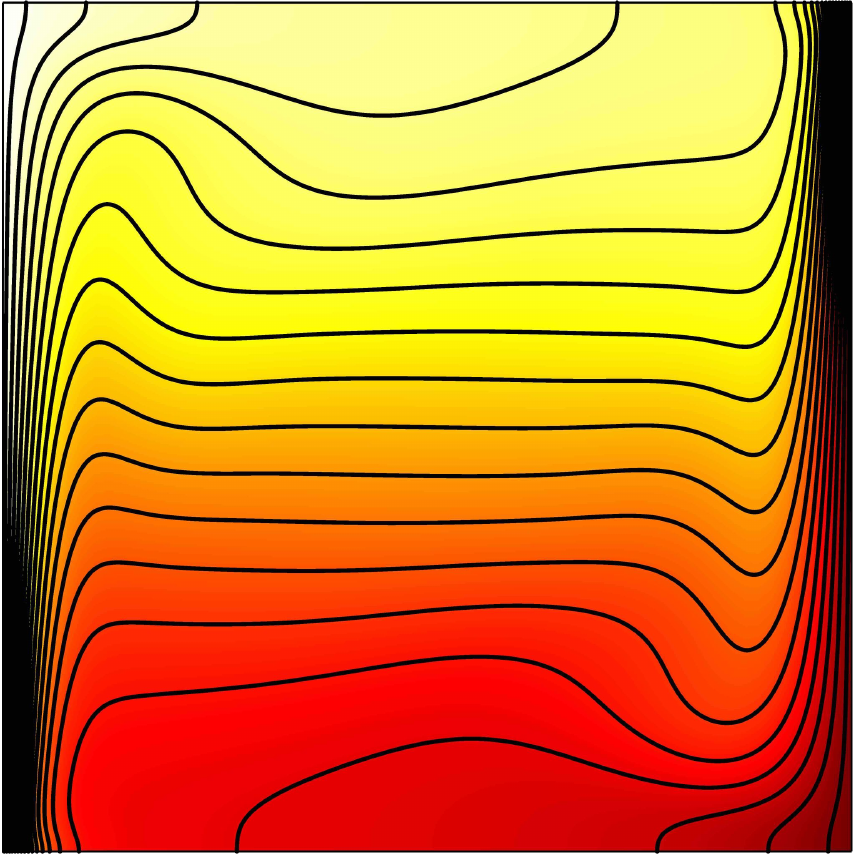}} \goodgap
\subfigure[$\Ray=\num{e8}$; $\tend=750$]{\includegraphics[width=0.285\textwidth]{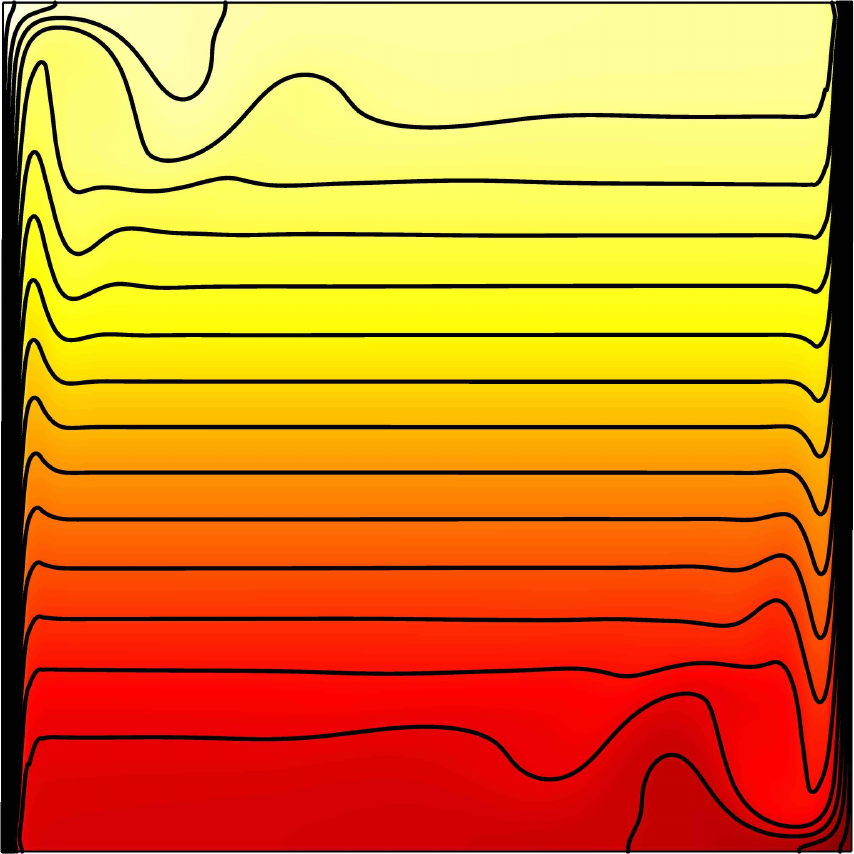}}
\caption{Temperature and isotherms of DHC simulation for $\Ray\in\left\{\num{e4},\num{e6},\num{e8}\right\}$. Results are computed with the $\rb{\Quad{-2}\slash\Quad{-1}}\wedge\Quad{-2}$ dG-FEM on $64\times 64$ meshes with stabilisation parameters $\gamma=\num{e5}$ and $\lambda=\num{e3}$.}
\label{fig:3}
\end{figure}

In Figure \ref{fig:3} the temperature fields and isotherms obtained by the same dG-FEM for different Rayleigh numbers are shown. Again, excellent agreement of our simulations with the reference solutions can be observed. As $\Ray$ increases the dominant mode for the heat transfer process changes from heat diffusion to heat convection. For $\Ray=\num{e4}$ the isotherms are more vertical than horizontal; especially near the vertical walls. For $\Ray=\num{e8}$, however, the isotherms form a strong thermal boundary layer at the vertical walls whilst being nearly horizontal in the rest of the domain. During the transition from diffusion-dominated to convection-dominated heat transfer, the isotherms for $\Ray=\num{e6}$ exhibit features of both modes. It is again important to note that although we used a coarse fixed-mesh without refinement, the sharp thermal and velocity boundary layers at the left and right walls are properly resolved by the proposed method.

\subsection{Benchmark quantities and comparative numerical schemes} 

In addition to verifying qualitatively that our numerical solution is in agreement with previous research we also want to regard the Nusselt number, representing heat flux across the cavity, as a quantity measuring the quality of the approximate temperature. It is important for the reliable application of numerical methods to consider such benchmark values and we note that in the context of conforming FEM \cite{GalvinEtAl12,JenkinsEtAl14}, the advantage of heavy grad-div stabilisation has only been demonstrated qualitatively without any quantitative analyses. \\

Therefore, consider the following Nusselt numbers, averaged over the mid-vertical line $x_1=0.5$ and averaged over the whole domain, which coincide for the true solution of the DHC problem \cite{GjesdalEtAl06}:
\begin{subequations}
\begin{align}
 	\Nu\rb{0.5} &= \frac{1}{\alpha}\oint_0^1 \sqb{u_{1h}T_h-\alpha \frac{\partial T_h}{\partial x_1}}\rb{0.5,x_2}\,\drm x_2 \\
	\NuAvg &= \frac{1}{\alpha}\int_\Omega \sqb{u_{1h}T_h-\alpha \frac{\partial T_h}{\partial x_1}}\,\drm \x	
\end{align}	
\end{subequations}
Note that, for the sake of notation, we agree to consider the Nusselt numbers only for the time-independent solution and thereby omit an explicit dependence on $t$. Moreover, in order to obtain a measure for the quality of the approximation of the velocity field, we also compute the stream function $\psi$ from \cite{Quere91}
\begin{subequations} \label{eq:Poisson}
	\begin{empheq}[left=\empheqlbrace]{alignat=2} 
	-\Delta \psi &=\rb{\frac{\partial u_{2h}}{\partial x_1}-\frac{\partial u_{1h}}{\partial x_2}}
	\qquad && \text{in}~\Omega, \\
	\psi&=0 && \text{on}~\partial\Omega,
	\end{empheq} 
\end{subequations}
where the occurring velocities are again considered to be the stationary solutions of the DHC computations. Problem \eqref{eq:Poisson} is always solved by a standard conforming, fifth-order FEM on the same mesh as the corresponding DHC simulation. In particular, we analyse the maximum absolute value over the whole cavity of the stream function and the absolute value measured exactly in the middle, i.e.
\begin{align}
\abs{\psi}_\maxrm =\max_{\x\in\Omega}\abs{\psi\rb{\x}}\quad\text{and}\quad
\abs{\psi}_\midrm =\abs{\psi\rb{0.5,0.5}}.
\end{align}
We remark that, on the discrete level, $\abs{\psi}_\midrm$ does not necessarily have to coincide with the maximum value of the stream function $\abs{\psi}_\maxrm$. Lastly, we want to measure how good and in which sense the discrete velocity $\uu_h$ satisfies the divergence constraint \eqref{eq:EPConti} and thus assess the mass conservation properties of the particular method. First of all, we note that due to the only approximate continuity across inter-element boundaries in the discrete velocity space $\UU_h^k$, the considered dG-FEM methods are non-conforming w.r.t.\ the space
\begin{equation}
\Hdiv=
\left\{\vv\in\sqb{\Lp^2\rb{\Omega}}^2\colon\nabla\cdot\vv\in\Lp^2\rb{\Omega}\right\}	,
\end{equation}
since the jump of the normal component across interfaces of the discrete velocity generally does not vanish \cite{BastianRiviere03}. Thus, the resulting discrete velocity from a classical dG-FEM does not have a well-defined divergence in the whole space $\Lp^2\rb{\Omega}$ but only locally on each $K\in\T$ and we can only consider the broken $\Lp^2$-norm 
\begin{equation}
\norm{\nabla_h\cdot\uu_h}_0	=
\rb{\sum_{K\in\T}\int_K \abs{\nabla\cdot\uu_h}^2\,\drm\x}^{\nicefrac{1}{2}}.
\end{equation}
This quantity will be used as a measure for the conservation of mass for all methods which are only approximately $\Hdiv$-conforming. \\

In order to put the results from our proposed stabilised dG-FEM into perspective, we want to compare them with solutions obtained by different finite element schemes. The easiest comparison is to compute an approximate solution with a standard $\Hk^1{\rb{\Omega}}$-conforming FEM. To this end we use a classical Taylor--Hood method with globally continuous $\rb{\Quad{2}\slash\Quad{1}}\wedge\Quad{2}$ elements with strongly imposed boundary conditions. Note that such a method naturally is also $\Hdiv$-conforming and thus the divergence of the resulting discrete velocity field belongs to $\Lp^2\rb{\Omega}$ globally. A well-known major drawback is that this method suffers from poor mass conservation since the divergence constraint is fulfilled only in a weak sense \cite{GalvinEtAl12,Linke09}.\\

Furthermore, we want to consider a method which, from the discrete function spaces used, can be described as being located between a conforming FEM and the proposed fully non-conforming dG-FEM of Section \ref{sec:Numerics}. Originating in the work about LDG methods in \cite{CockburnEtAl05b} we consider the possibility of constructing an exactly divergence-free, inf-sup stable and $\Hdiv$-conforming dG-FEM as described in \cite{CockburnEtAl07,WangEtAl09}. A simple way to achieve this within our previously defined framework is to consider the semi-discrete variational formulation \eqref{eq:DGVarForm} with $\kpr=k$ and without grad-div and pressure stabilisation ($\lambda=\gamma=0$) and to adapt the corresponding approximation space for the velocity. Therefore, we use the Raviart--Thomas space $\RTquad{k}$ on quadrilateral elements \cite{BoffiEtAl13} and define the following discrete space for the velocity: 
\begin{equation}
\UU_\RT^k=
\left\{\vv_h\in\Hdiv\colon\restr{\vv_h}{K}\in\RTquad{k}\rb{K},~\forall K\in\T;~\restr{\vv_h\cdot\n}{\partial\Omega}=0\right\}	
\end{equation}
In such a way we guarantee perfect mass conservation since $\nabla\cdot~\UU_\RT^k\subseteq\PP_h^k$. Note that the strong incorporation of the no-penetration boundary condition $\vv_h\cdot\n=0$ in the velocity space $\UU_\RT^k$ is necessary for obtaining an exactly divergence-free solution \cite{CockburnEtAl07}. However, the tangential component of the velocity on the boundary is still only imposed in the weak sense as described in Section \ref{sec:Numerics}. We refer to this method as $\rb{\RTquad{k}\slash\Quad{-k}}\wedge\Quad{-k}$. Lastly, we want to emphasise that due to the exact fulfilment of the divergence constraint the velocity error for this method does not depend on the pressure error and thereby yields a pressure-robust method in the sense of \cite{JohnEtAl16,Linke14,LinkeEtAl16}.

\subsection{Assessment of differently stabilised dG-FEM} 

Unfortunately, a typical situation encountered in applications is that the underlying computational mesh is not sufficiently resolved in critical flow areas. Our aim is to present results for such sub-optimal cases. Therefore, three different equidistant meshes with $16\times 16$, $32\times 32$ and $64\times 64$ quadratic elements will be used for the dG-FEM simulations of the DHC problem. In regard to the solutions in Figures \ref{fig:2} and \ref{fig:3} for higher Rayleigh numbers we clearly observe, and fully intend, that the occurring velocity and thermal boundary layers are locally under-resolved. We believe that in such a situation the advantages of the proposed stabilised dG-FEM are surfacing more clearly. \\

Regarding Table \ref{tab:1} we see all the above introduced benchmark quantities for the DHC problem with $\Ray=\num{e8}$ computed for different numerical schemes. Here, the bold numbers correspond to the particular values which are closest to the reference values given in \cite{Quere91,GjesdalEtAl06}. Thereby, for each mesh, they indicate the method which performs best in terms of mass conservation, Nusselt numbers and stream function maxima, respectively. The reason for beginning with the highest Rayleigh number is motivated by this being the computationally most demanding situation. Considering the first block in Table \ref{tab:1}, Taylor--Hood (TH) and equal-order (EO) type dG-FEM with and without pressure jump and grad-div stabilisation are compared on a $32\times 32$ mesh. First of all, we observe that even on such a coarse mesh all solutions are quite acceptable for all considered dG-FEM. Furthermore, in terms of mass conservation, the same comparatively high amount of grad-div stabilisation ($\gamma=\num{e5}$) improves the TH type methods more efficiently than the EO type method. The $\Lp^2$-norm of the broken divergence of the former is reduced about one order of magnitude more, even though without grad-div stabilisation they are on a comparable level. Note that for the heavily pressure stabilised ($\lambda=\num{e3}$) TH type dG-FEM, computations with $\gamma=0$ did not converge at all\textemdash{}an observation for which we do not have an explanation. \\

Now, we regard the second block of Table \ref{tab:1} where different dG-FEM are compared on a finer $64\times 64$ mesh. Again, all methods perform well for the considered problem and the benchmark quantities are close to the reference values. Moreover, it can be observed that the EO order methods generally seem to yield  better results for the temperature but worse for the velocity discretisation. In view of the Nusselt numbers it is remarkable that heavy grad-div stabilisation improves the accuracy for the temperature approximation for all considered methods. For the velocity approximation, however, we see that for the TH method without pressure stabilisation and the EO method, a high $\gamma\gg 1$ seems to introduce too much dissipation, thereby decreasing the accuracy of the stream function values. By contrast, the pressure stabilised TH dG-FEM has the nice behaviour that a high grad-div stabilisation simultaneously improves both the mass conservation properties and the overall accuracy in terms of the benchmark quantities. For all methods results for a different maximum amount of grad-div stabilisation are shown, corresponding to the experimentally obtained maximum value for which simulations converged. Thus, the pressure stabilised TH type dG-FEM turns out to be the most robust method w.r.t.\ heavy grad-div stabilisation. However, we note that beyond a certain value, increasing the grad-div parameter only yields a smaller broken divergence but the other benchmark quantities basically remain the same. For example, this value is $\gamma=\num{e2}$ for the EO type dG-FEM.  \\

\begin{table}[h]
\caption{Benchmark quantities for the DHC problem with $\Ray=\num{e8}$. DG-FEM computations of Taylor--Hood (TH) and equal-order (EO) type on different meshes with different grad-div ($\gamma$) and pressure ($\lambda$) stabilisation are compared with the reference solutions given in \cite{Quere91,GjesdalEtAl06}. A FEM solution with $\rb{\Quad{2}\slash\Quad{1}}\wedge\Quad{2}$ elements without grad-div stabilisation is computed on a $100\times 100$ mesh and a $\Hdiv$-conforming solution with $\rb{\RTquad{2}\slash\Quad{-2}}\wedge\Quad{-2}$ elements is computed on a $60\times 60$ mesh. The bold numbers correspond to the particular values which are closest to the reference values.}
\label{tab:1}
\centering
\begin{tabular}{ccccccccc}
\toprule
 Mesh  & Type &$\lambda$ &$\gamma$	&$\norm{\nabla_h\cdot \uu_h}_0$ & $\Nu\rb{0.5}$	& $\NuAvg$	&$\sqrt{\Ray}\abs{\psi}_\maxrm$ & $\sqrt{\Ray}\abs{\psi}_\midrm$ 					\\
\otoprule
$32\times 32$
 &TH &\num{0}   	&\num{0}		&\num{0.0764}			&28.161 			&28.195		 &54.756 		&{\bf 52.432}\\
 &TH &\num{0}   	&\num{e5}	&\num{1.99e-9}			&28.353 			&28.356 	 &54.723 	   	&52.467\\
 &TH &$\num{e3}$ 	&\num{0.01}	&\num{0.0146}			&28.118 			&28.145 	 &54.717 	   	&52.459\\
 &TH &$\num{e3}$ 	&\num{e5}	&{\bf\num{1.88e-9}}		&28.182 			&28.184 	 &{\bf54.636}	&52.443\\
 &EO &$\num{1}$ 	&\num{0}		&\num{0.0625}			&28.541 			&{\bf28.545} &55.041 	   	&52.746\\
 &EO &$\num{1}$ 	&\num{e5}	&\num{1.82e-8}			&{\bf28.639}		&28.338 	 &55.530 	   	&52.068\\
\midrule
$64\times 64$
 &TH &\num{0}		&\num{0}		&\num{0.0290}			&30.061 			&30.063   	 &53.835 	   	&52.206 \\
 &TH &\num{0}		&\num{e5}	&\num{9.77e-10}			&30.071 			&30.069 	 &53.798 	   	&52.165 \\
 &TH &$\num{e3}$	&\num{0}		&\num{0.0560}			&29.931 			&29.944  	 &53.635 	   	&52.020 \\
 &TH &$\num{e3}$	&\num{e5}	&\num{8.40e-10}			&30.064 			&30.062 	 &{\bf53.847}	&52.217 \\
 &TH &$\num{e3}$	&\num{e9}	&{\bf\num{8.40e-14}}		&30.066 			&30.062 	 &{\bf53.847}	&{\bf52.220}\\
 &EO &$\num{1}$		&\num{0}		&\num{0.0306}			&30.098 			&{\bf30.099} &53.800 	   	&52.168\\
 &EO &$\num{1}$		&\num{e2}	&\num{2.13e-5}			&30.137 		 	&30.029 	 &53.800 	   	&51.782\\
 &EO &$\num{1}$		&\num{e7}	&\num{2.15e-10}			&{\bf30.138}		&30.029 	 &53.800 	   	&51.779 \\
\midrule
 	&&& FEM			&\num{0.0272}		&30.174			&30.174		 	&53.785		 &52.255	\\	
 	&&& $\Hdiv$		&\num{3.43e-15}		&30.062			&30.063	 	 	&53.794		 &52.155	\\
 	&&& \cite{Quere91} 		   		&-					&30.225  		&30.225		 &53.85  		&52.32	\\
	&&& \cite{GjesdalEtAl06}   		&-					&30.223  		&30.223 	 &53.84  		&52.32 	\\
\bottomrule
\end{tabular}
\end{table}

Regarding the FEM without grad-div stabilisation and the $\Hdiv$-conforming method in the last block of Table \ref{tab:1}, we first note that using $\rb{\Quad{2}\slash\Quad{1}}\wedge\Quad{2}$ elements on a $100\times 100$ mesh and $\rb{\RTquad{2}\slash\Quad{-2}}\wedge\Quad{-2}$ elements on a $60\times 60$ mesh yields a comparable number of degrees of freedom (DOF) as the considered dG-FEM on a $64\times 64$ mesh. Furthermore, we directly observe that on such a fine mesh the standard FEM yields a slightly better temperature approximation as our proposed dG-FEM but both the $\Lp^2$-norm of the divergence and the stream function values from the velocity approximation are considerably worse. In the next subsection, the impact of additional grad-div stabilisation for this method is investigated separately. The $\Hdiv$-conforming and exactly divergence-free method, on the other hand, by construction yields perfect mass conservation and also the Nusselt numbers and stream function values agree very well with both the reference values and our stabilised dG-FEM.\\

However, the  major drawback of this method is much more subtle. Even though we have about the same number of DOF as for the proposed TH and EO dG-FEM, the $\rb{\RTquad{2}\slash\Quad{-2}}\wedge\Quad{-2}$ method takes more than twice the computing time to finish the DHC simulation. We believe that the reason for this disadvantage is based on the fact that LDG methods have a larger stencil than interior penalty methods and thus are generally up to 2.5 times less efficient from a computational point of view \cite{Castillo02}.\\

Furthermore, it is remarkable that the heavily grad-div stabilised ($\gamma=\num{e9}$) TH dG-FEM yields a broken divergence whose magnitude is comparable to this exactly divergence-free method. As the results from the $\Hdiv$-conforming and our proposed stabilised dG-FEM are qualitatively comparable we infer that the  advantage of pressure-robustness is not significant for the considered DHC problem with a high Rayleigh number. From this comparison we deduce that the proposed stabilised dG-FEM are indeed very well-suited methods for solving weakly non-isothermal computational fluid dynamics problems with boundary layers on under-resolved and non-adapted meshes in the sense of both overall accuracy and efficiency. \\

\begin{figure}[h]
\centering
	\subfigure[TH, $\lambda=\num{e3}$, $\gamma=0$ \label{fig:4a}]{\includegraphics[height=0.165\textheight]{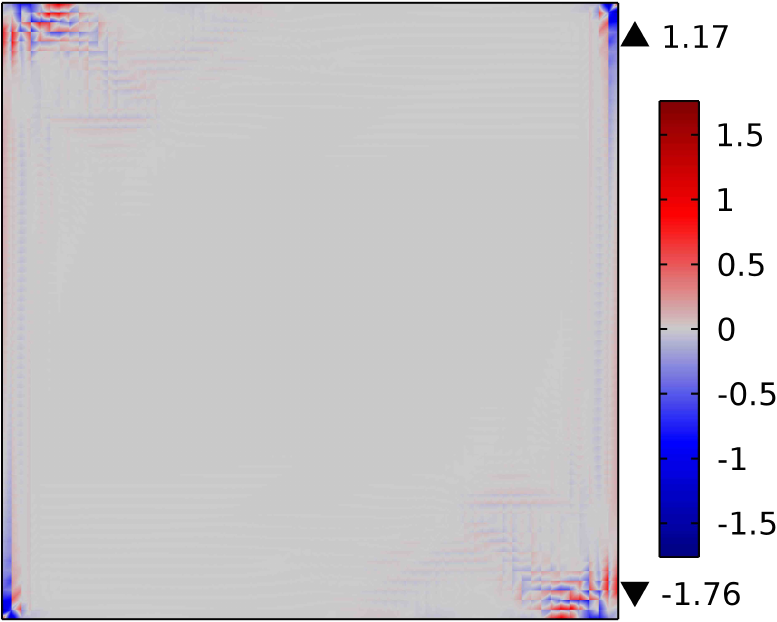}}	\goodgap
	\subfigure[TH, $\lambda=\num{e3}$, $\gamma=\num{e9}$ \label{fig:4b}]{\includegraphics[height=0.165\textheight]{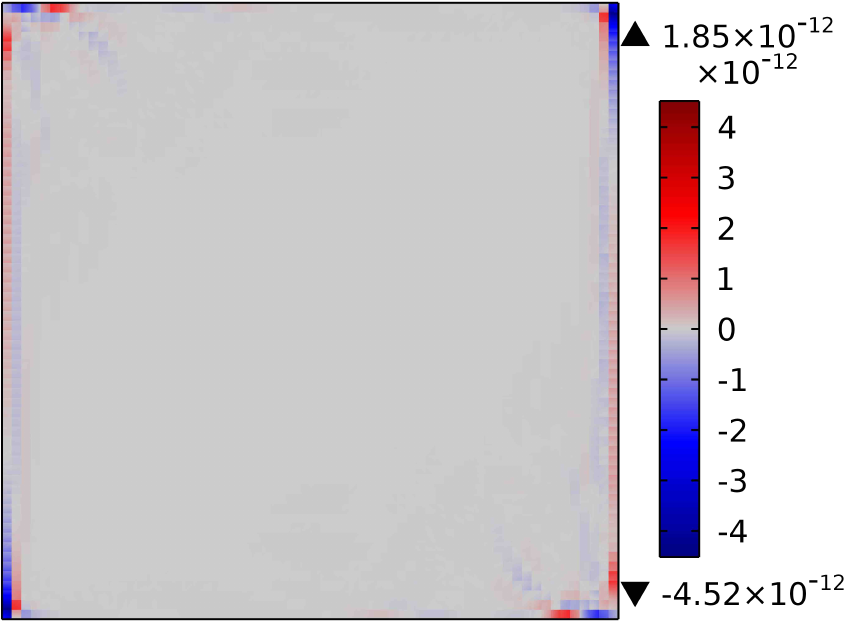}}	\goodgap
	\subfigure[EO, $\lambda=\num{1}$, $\gamma=\num{e7}$ \label{fig:4c}]{\includegraphics[height=0.165\textheight]{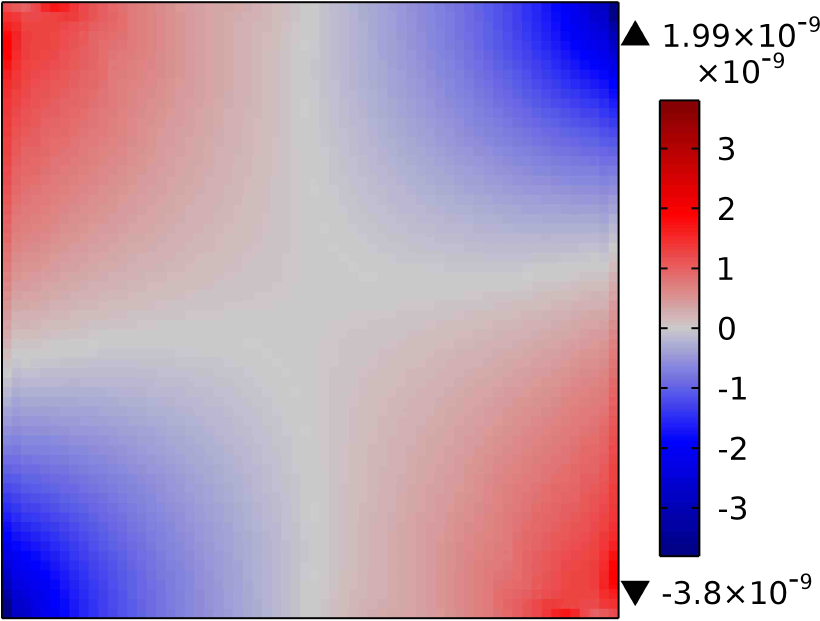}}
\caption{Pointwise broken divergence $\nabla_h\cdot\uu_h$ resulting from differently stabilised dG-FEM on the $64\times 64$ mesh for $\Ray=\num{e8}$. Note the different scales of the legends. The upward and downward pointing black triangles indicate the maximum and minimum value attained over the whole domain, respectively.}
\label{fig:4}
\end{figure}

The last comparison gives a recommendation between TH and EO type stabilised dG-FEM. We already saw that the temperature approximation of EO type dG-FEM generally is better whereas the velocity approximation is worse compared to TH type dG-FEM. Also, in terms of the $\Lp^2$-norm of the broken divergence, especially pressure stabilised TH dG-FEM with heavy grad-div stabilisation are superior compared to EO dG-FEM, which also generally yield a higher number of DOF for the same problem. Additionally, we regard the pointwise broken divergence $\nabla_h\cdot\uu_h\rb{\x}$ for different methods in Figure \ref{fig:4}. The divergence in Figure \ref{fig:4a} results from a pressure stabilised TH dG-FEM without grad-div stabilisation and we see that poor mass conservation can be observed primarily in the boundary layer and towards the corners of the cavity. We note that EO type dG-FEM without grad-div stabilisation yield a very similar broken divergence which is not shown. Regarding Figure \ref{fig:4b} the general pattern of the broken divergence remains unaffected if heavy grad-div stabilisation ($\gamma=\num{e9}$) is used for the TH dG-FEM. The magnitude of the broken divergence, however, is reduced by twelve orders of magnitude, thereby yielding an essentially pointwise divergence-free velocity approximation. If heavy grad-div stabilisation ($\gamma=\num{e7}$) is used for the EO type dG-FEM the magnitude of the broken divergence can also be reduced significantly; cf. Figure \ref{fig:4c}. However, the regions with non-zero divergence are not confined in the boundary layer anymore which, physically speaking, is a significant drawback of the EO type dG-FEM. Indeed, blue regions with negative divergence and red regions with positive divergence correspond to mass sinks and sources, respectively \cite{Durst08}. Furthermore, it has not been possible to use a grad-div parameter higher than $\gamma=\num{e7}$ for the EO type methods. \\

Due to all the above mentioned reasons, we rate the pressure stabilised TH dG-FEM as the most promising method for solving thermally-coupled incompressible flow problems. Therefore, exclusively a TH type dG-FEM with stabilisation parameters $\gamma=\num{e5}$ and $\lambda=\num{e3}$ is used for the computation of the benchmark quantities in Table \ref{tab:2} for smaller Rayleigh numbers. \\

Regarding the computed benchmark quantities in Table \ref{tab:2} we observe an excellent agreement with the reference values given in \cite{Davis83,Quere91,GjesdalEtAl06,SaitohHirose89} even on coarse meshes. Additionally, due to the heavy grad-div stabilisation ($\gamma=\num{e5}$) the broken divergence is reduced significantly. The reason for showing these results is to illustrate that the proposed stabilised dG-FEM performs well for diffusion-dominated problems. \\

\begin{table}[h]
\caption{Benchmark quantities for the DHC problem with $\Ray\in\{\num{e4},\num{e6}\}$. Our own results, obtained by the Taylor--Hood type dG-FEM with $\gamma=\num{e5}$ and $\lambda=\num{e3}$ on different meshes, are compared with the reference solutions given in \cite{Davis83,Quere91,GjesdalEtAl06,SaitohHirose89}.}
\label{tab:2}
\centering
\begin{tabular}{ccccccc}
\toprule
$\Ray$	& Mesh &$\norm{\nabla_h\cdot \uu_h}_0$	  	& $\Nu\rb{0.5}$	& $\NuAvg$	&$\sqrt{\Ray}\abs{\psi}_\maxrm$ & $\sqrt{\Ray}\abs{\psi}_\midrm$ 				\\
\otoprule
\num{e4}
&$16\times 16$	  	&\num{2.26e-8}	&2.2384  &2.2433  &5.0747 &5.0747 	\\
&$32\times 32$	  	&\num{6.05e-9}	&2.2435  &2.2447  &5.0742 &5.0742 	\\
&$64\times 64$	  	&\num{1.54e-9}	&2.2445  &2.2448  &5.0738 &5.0738 	\\
\cmidrule(l){2-7} 
	& \cite{Davis83}: 		  &-		&\num{2.243}  &\num{2.243}  &-	&\num{5.071} \\
	& \cite{SaitohHirose89}:  &-		&\num{2.2430} &\num{2.2424} &-	&\num{5.0731}\\
\midrule
\num{e6}
&$16\times 16$	  	&\num{1.36e-8}	&8.6982  &8.6892  &16.830 &16.359 \\
&$32\times 32$	  	&\num{5.18e-9}	&8.8121  &8.8108  &16.816 &16.382 \\
&$64\times 64$	  	&\num{1.47e-9}	&8.8241  &8.8239  &16.815 &16.388 \\
\cmidrule(l){2-7} 	
	& \cite{Quere91}: 		&-	&\num{8.825}  &\num{8.825}	&\num{16.811}	&\num{16.386} \\
	& \cite{GjesdalEtAl06}: &-  	&\num{8.825}  &\num{8.825}	&\num{16.809}	&\num{16.384} \\
\bottomrule
\end{tabular}
\end{table}

\subsection{Mass conservation for conforming FEM} 

In this last subsection on the DHC problem we want to investigate the impact of heavy grad-div stabilisation on standard, conforming FEM as suggested by \cite{GalvinEtAl12,JenkinsEtAl14}. Therefore, we computed benchmark quantities for varying $\gamma$ with a standard conforming FEM with $\rb{\Quad{2}\slash\Quad{1}}\wedge\Quad{2}$ elements on a $50\times 50$ mesh for the DHC problem with $\Ray=\num{e8}$. This method yields a number of DOF which is comparable to the dG-FEM on the $32\times 32$ mesh. For the previously proposed stabilised dG-FEM, we saw that heavy grad-div stabilisation improves both the mass conservation and the overall accuracy of the method. As it turns out, the simultaneity of improving these two properties is by no means self-evident. \\

Indeed, regarding Table \ref{tab:3} we see that the FEM without grad-div stabilisation yields acceptable benchmark quantities which are in fact comparable to the ones obtained by the corresponding dG-FEM. Additionally, we observe the (not surprising but still positive) fact that heavy grad-div stabilisation also yields approximate velocity fields with a significantly decreased divergence. However, the accuracy of both temperature and velocity approximation apparently deteriorates completely for large values of $\gamma$. Therefore, we infer that for standard FEM, simultaneously having good mass conservation properties and a good accuracy is not possible through heavy grad-div stabilisation. \\

\begin{table}[h]
\caption{Impact of grad-div stabilisation on a standard conforming FEM with $\rb{\Quad{2}\slash\Quad{1}}\wedge\Quad{2}$ elements on a $50\times 50$ mesh for the DHC problem with $\Ray=\num{e8}$.  For comparison, the reference values from \cite{Quere91,GjesdalEtAl06} are shown. The bold numbers correspond to the particular values which are closest to the reference values.}
\label{tab:3}
\centering
\begin{tabular}{cccccc}
\toprule
 $\gamma$   	&$\norm{\nabla\cdot \uu_h}_0$ & $\Nu\rb{0.5}$	& $\NuAvg$	&$\sqrt{\Ray}\abs{\psi}_\maxrm$ & $\sqrt{\Ray}\abs{\psi}_\midrm$ 					\\
\otoprule
$\num{0}$	 	&0.1024					&{\bf 29.332}	&{\bf 29.364} 	&{\bf 52.794}	&{\bf 51.241} 	\\
$\num{0.25}$	&0.0183					&28.502 	 		&28.496			&51.332 		&49.691 		\\
$\num{1}$	 	&0.0114					&27.058 	 		&27.053			&50.752 		&48.904 		\\
$\num{2}$	 	&\num{7.80e-3}			&26.218 	 		&26.215			&50.909 		&48.889 		\\
$\num{5}$	 	&\num{4.01e-3}			&25.344 	 		&25.342			&51.301 		&49.105 		\\
$\num{e2}$	 	&\num{2.45e-4} 			&24.491 	 		&24.491			&51.853 		&49.494 		\\
$\num{e5}$	 	&{\bf \num{2.48e-7}}	 	&24.436 	 		&24.436			&51.892 		&49.524 		\\
\midrule 	
	 \cite{Quere91}: 		&-	&\num{30.225}  &30.225		 &\num{53.85}  &\num{52.32} \\
	 \cite{GjesdalEtAl06}: &-	&\num{30.223}  &\num{30.223} &\num{53.84}  &\num{52.32} \\
\bottomrule
\end{tabular}
\end{table}

\subsection{Personal recommendation} 

Let us briefly summarise the results of the DHC simulations and comment on the suitability of the proposed stabilised methods. Referring to Table 1, we saw that the EO type dG-FEM in general yields a more accurate temperature approximation whereas the grad-div and pressure stabilised TH type dG-FEM is superior in terms of the velocity approximation. Therefore, whenever the focus of the particular simulation is on temperature-related phenomena, we recommend to use the EO type dG-FEM with mild grad-div stabilisation (e.g.\ $\gamma=\num{E2}$). With this choice, the compromise to make lies in the slightly worse velocity approximation and the larger number of DOF due to the richer pressure FE space. However, we favour the TH type dG-FEM with both pressure and grad-div stabilisation due to its robustness with respect to large stabilisation parameters. Indeed, in our opinion, using the TH type dG-FEM with stabilisation parameters $\gamma=\num{e5}$ and $\lambda=\num{e3}$ is a safe and robust choice for a wide variety of non-isothermal fluid flow problems which yields an accurate approximation with excellent mass conservation properties and comparably few DOF. Therefore, this method is chosen exclusively for the next section where we consider a phase change problem with moving interior layers.

\section{Moving interior layers: Two-phase model for melting of pure gallium} \label{sec:Gallium}

For the purpose of further assessing the quality and performance of the proposed stabilised dG-FEM, the melting of pure gallium in a differentially heated enclosure is considered. This problem has been used frequently for the assessment of numerical schemes involving melting and solidification processes with a moving interior layer. For solid/liquid phase change processes in general and for gallium melting in particular, there is only few experimental data available \cite{GauViskanta86,CampbellKoster94}. Therefore, our results are compared to the ones obtained by \cite{BelhamadiaEtAl12,CagnoneEtAl14,StellaGiangi00,HannounEtAl03} which are numerical results published over the last 15 years using finite element, finite volume and also discontinuous Galerkin methods. Additionally, a comprehensive analysis of different time-stepping schemes and step sizes can be found in \cite{EvansKnoll07} where it is shown that, in the framework of phase change problems, the usage of BDF(2) is clearly preferable to lower-order schemes as for example the implicit Euler method. In agreement with the previous section on the DHC problem, exclusively a TH type dG-FEM with stabilisation parameters $\gamma=\num{e5}$ and $\lambda=\num{e3}$ is used for the computation of the gallium problem.

\subsection{Problem statement}

Suppose we have a block of solid gallium with melting point $T_f$, initially held at a constant temperature $T_0<T_f$ in a rectangular cavity of width $W$ and height $H$. The top and bottom walls of this cavity are assumed to be adiabatic. Then, start to increase the temperature at the left wall to $T_\hot >T_f$ whilst maintaining the temperature on the right wall at $T_\cold=T_0<T_f$. The hot wall, providing a temperature above the melting point, causes the gallium to melt and form a liquid phase across the left wall. For the fluid phase the no-slip condition is imposed on all walls of the cavity. Furthermore assume that gravity, inducing a motion in the melt flow, acts in the negative $x_2$-direction. A suitable mathematical model describing this problem is given by the enthalpy-porosity method \eqref{eq:EnthalpyPorosityEQs}. All relevant physical properties of gallium together with the other system parameters can be found in Table \ref{tab:4}. \\
 
\begin{table}[h]
\caption{Physical properties and system parameters for the melting of pure gallium.}
\label{tab:4}
\centering
\begin{tabular}{cccc}
\toprule
Property/Parameter			& Symbol						& Value				& Unit						\\
\otoprule
Density  					& $\rho$					& \num{6093}			& $\si{kg.m^{-3}}$			\\
Kinematic viscosity 		& $\nu$						& \num{2.97e-7}		& $\si{m^2.s^{-1}}$			\\
Thermal expansion			& $\beta$					& \num{1.2e-4}		& $\si{K^{-1}}$				\\
Gravitational acceleration	& $g$						& 10					& $\si{m.s^{-2}}$			\\
Specific heat capacity 		& $c_p$						& 381.5				& $\si{J.kg^{-1}.K^{-1}}$	\\
Thermal conductivity 		& $\kappa$					& 32					& $\si{W.m^{-1}.K^{-1}}$	\\
Latent heat of fusion		& $L_f$						& \num{80160}		& $\si{J.kg^{-1}}$			\\
Temperature of fusion		& $T_f$						& 302.78				& $\si{K}$					\\
\midrule
Velocity attenuation		& $C_0$						& \num{e8}			& $\si{kg.m^{-3}.s^{-1}}$	\\
Security parameter			& $b$						& \num{e-3}			& $\si{1}$					\\
Hot wall temperature		& $T_\hot$					& 311				& $\si{K}$					\\
Cold wall temperature 		& $T_\cold$					& 301.3				& $\si{K}$					\\
Reference temperature 		& $T_\refrm$				& 301.3				& $\si{K}$					\\
\midrule
Reference temperature difference & $\Delta {T_\refrm}$	& 9.7				& $\si{K}$					\\
Cavity width				& $W$						& \num{2e-2}			& $\si{m}$					\\
Cavity height				& $H$						& \num{6.35e-2}		& $\si{m}$					\\
Prandtl number				& $\Pra$					& 0.021658			& $\si{1}$					\\
Rayleigh number				& $\Ray$					& \num{7.2879e5}		& $\si{1}$					\\
\bottomrule
\end{tabular}
\end{table}

In order to minimise the computational costs, the cavity considered here is narrower compared to the one used for example in \cite{HannounEtAl03}. As a matter of fact, $W$ given in Table \ref{tab:4} corresponds to only $\SI{22.5}{\percent}$ of the width of the cavity in \cite{HannounEtAl03}. However, the focus in this work is on the early melting process and therefore it suffices to consider such a downsized domain. Note that in accordance with the literature, the reference length $L_\refrm$ used for computing the Grashof and thus the Rayleigh number is chosen to be the cavity height $H$. The closing initial and boundary conditions are specified in the following:
\begin{enumerate}[label=(\roman*)]
\item Initial values $u_1\equiv u_2\equiv p \equiv 0$ and $T=T_0=T_\cold$ on $\Omega=\rb{0,W}\times\rb{0,H}$ at $t=0$.
\item Dirichlet conditions $g_D^T=T=T_\hot$ on $x_1=0$ and $g_D^T=T=T_\cold$ on $x_1=W$ for all $0\leqslant x_2\leqslant H$ for the temperature. Define by $\Gamma_D=\left\{\rb{x_1,x_2}\in\Omega\colon x_1=0~\text{and}~x_1=W\right\}$ the Dirichlet part of the boundary $\partial\Omega$ with prescribed Dirichlet boundary condition $g_D^T$ on $\Gamma_D$.
\item Homogeneous Neumann conditions $g_N^T=\frac{\partial T}{\partial x_2}=0$ on $x_2=0$ and $x_2=H$ for all $0\leqslant x_1\leqslant W$ for the temperature. Define by $\Gamma_N=\left\{\rb{x_1,x_2}\in\Omega\colon x_2=0~\text{and}~x_2=H\right\}$ the Neumann part of the boundary $\partial\Omega$ with prescribed Neumann boundary condition $g_N^T$ on $\Gamma_N$.
\end{enumerate}
The space semi-discrete dG-FEM formulation is given by \eqref{eq:DGVarForm}. In accordance with Section \ref{sec:Numerics} all simulations are carried out on non-adapted meshes with quadratic elements and moreover, referring to \cite{EvansKnoll07}, in order to minimise the impact of the temporal discretisation we restrict the maximum time step of the BDF(2) solver to $\SI{0.05}{s}$. Again, to ensure the compatibility of boundary and initial conditions for the temperature, the hot wall temperature is ramped up smoothly from $T_0$ to $T_\hot$ during the first $\SI{0.1}{s}$ of simulation.

\subsection{Results and mesh convergence}
\begin{figure}[h]
\centering
\subfigure[$\SI{9}{s}$]{\includegraphics[height=0.4\textheight]{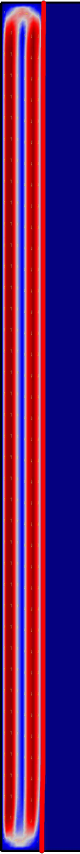}}		\goodgap 
\subfigure[$\SI{20}{s}$]{\includegraphics[height=0.4\textheight]{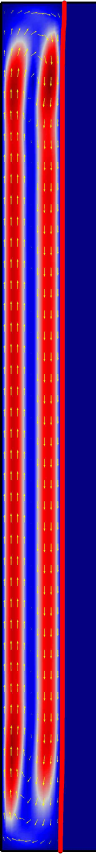}}		\goodgap 
\subfigure[$\SI{32}{s}$]{\includegraphics[height=0.4\textheight]{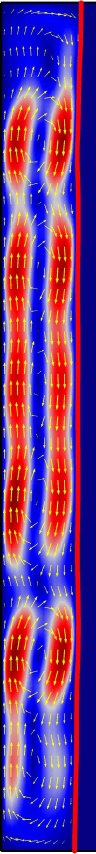}} 	\goodgap 
\subfigure[$\SI{36}{s}$]{\includegraphics[height=0.4\textheight]{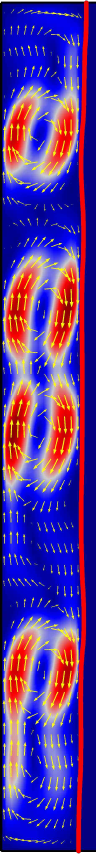}} 	\goodgap 
\subfigure[$\SI{42}{s}$]{\includegraphics[height=0.4\textheight]{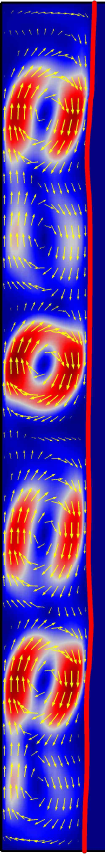}} 	\goodgap 
\subfigure[$\SI{50}{s}$]{\includegraphics[height=0.4\textheight]{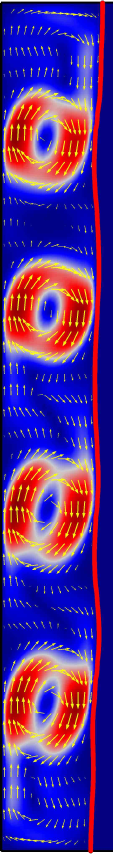}} 	\goodgap
\subfigure[$\SI{65}{s}$]{\includegraphics[height=0.4\textheight]{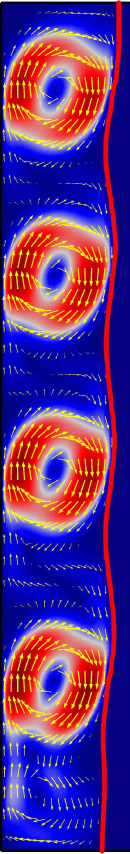}} 	\goodgap  
\subfigure[$\SI{85}{s}$]{\includegraphics[height=0.4\textheight]{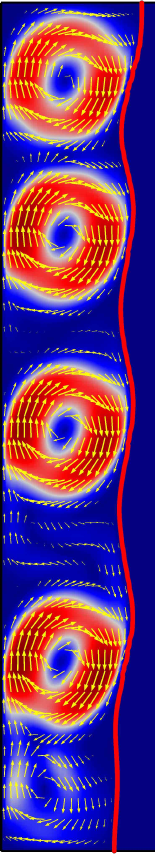}}	
\caption{Evolution of velocity magnitude and phase boundary ($0.5$-contour of $\phi$) for $\Delta T_f=\SI{0.25}{K}$ shown at different time instances. Computed by the TH dG-FEM with $\gamma=\num{e5}$ and $\lambda=\num{e3}$ on a $h=\SI{0.4}{mm}$ mesh. In each case only the left part of the enclosure, adjusted to the melt flow of the liquid phase, is shown.}
\label{fig:5}
\end{figure}

In Figure \ref{fig:5} one can see the velocity field and phase boundary (red line) at different time steps obtained by the TH type dG-FEM with $\gamma=\num{e5}$ and $\lambda=\num{e3}$ on a mesh with $h=\SI{0.4}{mm}$. For the melting range we choose $\Delta T_f=\SI{0.25}{K}$ which turns out to be an appropriate choice. However, in the next subsection we conduct a study for different melting ranges that corresponds to moving interior layers with different sharpnesses. \\

After the initial $t=\SI{9}{s}$ of heating from the left wall, the gallium develops a small liquid region along this wall. There is one big circulation where the fluid rises at the hot wall and drops at the phase boundary. At $t=\SI{20}{s}$ a slight tendency of the flow to develop two vortices at the top and bottom can be seen. Proceeding to $t=\SI{32}{s}$ we observe that the amount of liquid gallium increases and stand-alone vortices develop at both top and bottom of the enclosure where the melt accumulates. At $t=\SI{36}{s}$ these two vortices are visible more clearly and a tendency of the flow for developing two more vortices between them can be seen. As the evolution proceeds to $t=\SI{42}{s}$ all four vortices are amplified and are clearly separated from each other. At $t=\SI{50}{s}$ four concise vortices are developed in the enclosure whose size increases steadily towards $t=\SI{65}{s}$. Additionally to growing, buoyancy forces due to the Boussinesq term let them rise and at $t=\SI{85}{s}$ signs of another small vortex appear at the bottom of the enclosure. The velocity fields obtained by this proposed dG-FEM agree excellently with \cite{BelhamadiaEtAl12,CagnoneEtAl14,StellaGiangi00,HannounEtAl03} although the mesh for the present computation is considerably coarser. Indeed, our $h=\SI{0.4}{mm}$ mesh consists of \num{7950} quadratic elements whereas, for example in \cite{HannounEtAl03}, a central finite volume method with a fully implicit Euler scheme is used on a mesh with \num{113400} elements (interpolated value to comply with our smaller enclosure). This means that our mesh is more than 14 times coarser but nonetheless, with the proposed stabilised dG-FEM, still qualitatively yields the same results.\\

\begin{figure}[h]
\centering
\subfigure[Mesh 1]{\includegraphics[height=0.4\textheight]{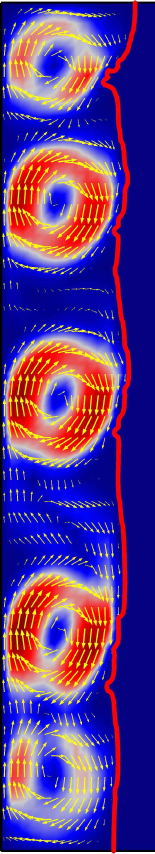}}	\goodgap 
\subfigure[Mesh 2]{\includegraphics[height=0.4\textheight]{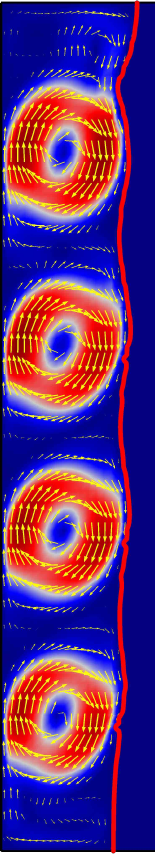}}	\goodgap 
\subfigure[Mesh 3]{\includegraphics[height=0.4\textheight]{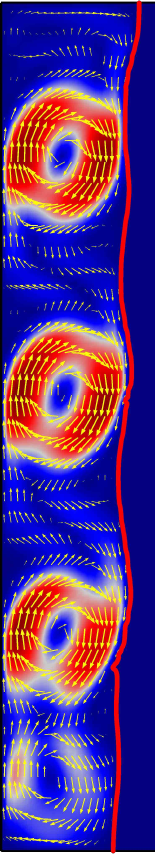}}	\goodgap 
\subfigure[Mesh 4]{\includegraphics[height=0.4\textheight]{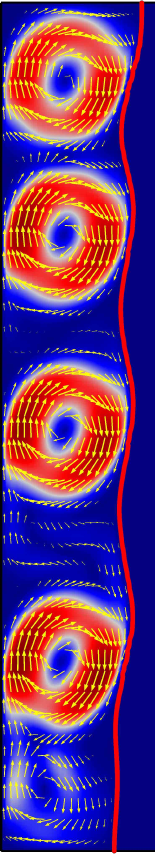}} 	\goodgap 
\subfigure[Mesh 5]{\includegraphics[height=0.4\textheight]{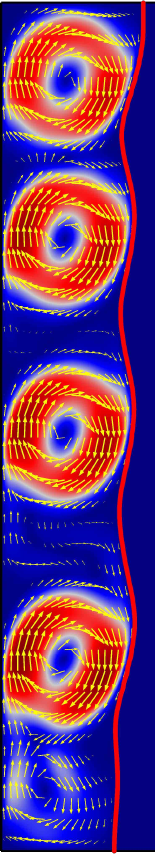}} 	\goodgap 
\subfigure[Mesh 6]{\includegraphics[height=0.4\textheight]{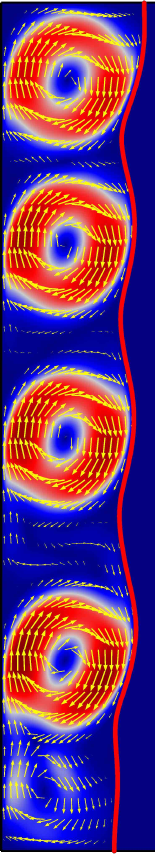}} 	
\caption{Mesh sensitivity analysis. Velocity magnitude and phase boundary ($0.5$-contour of $\phi$) at $t=\SI{85}{s}$ for $\Delta T_f=\SI{0.25}{K}$. Computed by the TH type dG-FEM with $\gamma=\num{e5}$ and $\lambda=\num{e3}$ on different meshes. In each case, only about the left half of the enclosure is shown. The meshes are defined in Table \ref{tab:5}.}
\label{fig:6}
\end{figure}

We briefly comment on similarities and differences between our proposed method and the one used in \cite{CagnoneEtAl14}. In \cite{CagnoneEtAl14} a SIP dG-FEM with BDF(2) time-stepping is considered for the enthalpy-porosity model, as well. However, in contrast to this work, the energy equation is written in terms of the enthalpy as primary variable. Concerning the numerical method, also a fully-coupled Taylor--Hood type dG-FEM is employed; but on an unstructured simplicial mesh which is refined towards the hot wall. Unfortunately, the order of the particular interpolation spaces is not mentioned. Furthermore, neither grad-div nor pressure stabilisation is considered. Regarding the simulation results, we note that the position of the vortices in \cite{CagnoneEtAl14} is slightly different from the position of the vortices in this work; see Figure \ref{fig:5}. However, our results coincide very well with the other references \cite{BelhamadiaEtAl12,StellaGiangi00,HannounEtAl03}. Therefore, we infer that our stabilised dG-FEM produces results which agree better with existing reference solutions for the problem of gallium melting. \\

The next step is to verify that the above shown results are robust against mesh refinement, thereby justifying the validity of the analysis. Therefore, we compare the velocity fields at $t=\SI{85}{s}$ obtained by the same dG-FEM on the six different meshes summarised in Table \ref{tab:5}. But first of all, we note that it is remarkable that the gallium simulation can be computed on the extremely coarse $\SI{0.8}{mm}$ mesh at all. This natural treatment of moving interior layers is clearly an advantage of stabilised dG-FEM. Regarding Figure \ref{fig:6} we observe that all meshes yield a flow structure with three to five vortices at different locations in the enclosure. However, only for $h\leqslant\SI{0.4}{mm}$ four vortices are basically fixed in space and thus independent of mesh refinement. This situation is in agreement with the literature and therefore, we infer that the gallium simulation with a melting range of $\Delta T_f=\SI{0.25}{K}$ is close enough to being mesh-converged such that all subsequent simulations are carried out on the $h=\SI{0.4}{mm}$ mesh. \\

\begin{table}[h]
\caption{Mesh size, number of mesh elements and number of DOF for the mesh sensitivity analysis of gallium melting.}
\label{tab:5}
\centering 
\begin{tabular}{ccccccc} 
\toprule
Mesh				&1 			 &2			  &3 			 &4 			   &5 			 &6				\\ 
\otoprule
$h$ $\sqb{\si{mm}}$	&\num{0.8}	 &\num{0.6}	  &\num{0.5}     &\num{0.4} 	   &\num{0.3}	 &\num{0.2}		\\ 
Mesh elements		&\num{2000}	 &\num{3604}  &\num{5080}	 &\num{7950}	   &\num{14204}  &\num{31800} 	\\
DOF					&\num{62001} &\num{111725} &\num{157481} &\num{246451} &	\num{440325} &\num{985801}	\\ 
\bottomrule
\end{tabular}
\end{table}

\subsection{Sharpness of interior layer}

\begin{figure}[h]
\centering
\subfigure[$\SI{2}{K}$]{\includegraphics[height=0.4\textheight]{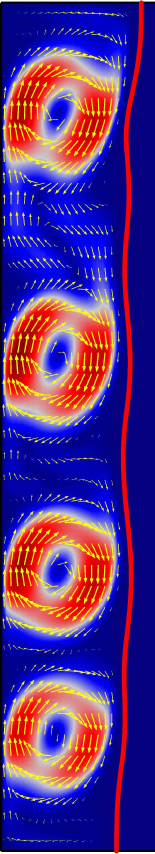}}		\goodgap 
\subfigure[$\SI{1}{K}$]{\includegraphics[height=0.4\textheight]{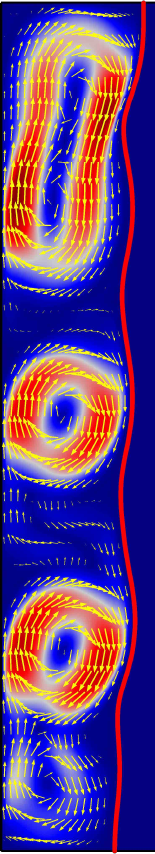}} 		\goodgap 
\subfigure[$\SI{0.8}{K}$]{\includegraphics[height=0.4\textheight]{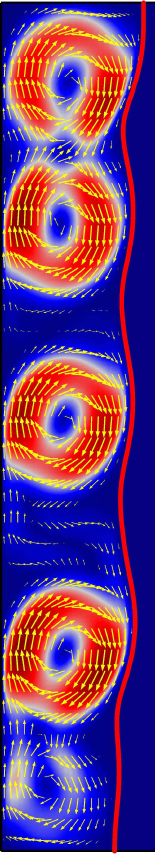}} 	\goodgap 
\subfigure[$\SI{0.25}{K}$]{\includegraphics[height=0.4\textheight]{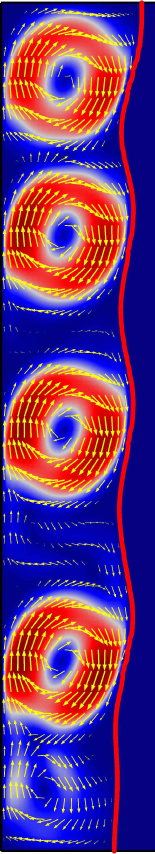}} 	\goodgap 
\subfigure[$\SI{0.175}{K}$]{\includegraphics[height=0.4\textheight]{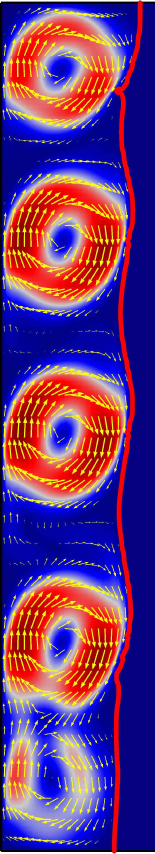}} \goodgap
\subfigure[$\SI{0.125}{K}$]{\includegraphics[height=0.4\textheight]{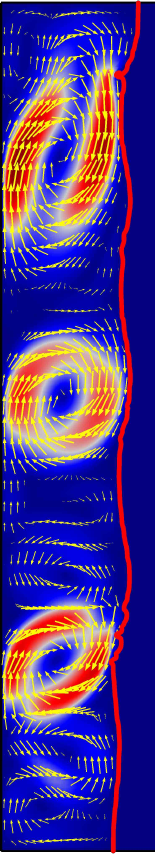}}
\caption{Velocity magnitude and phase boundary ($0.5$-contour of $\phi$) at $t=\SI{85}{s}$  computed on a $h=\SI{0.4}{mm}$ mesh by the TH type dG-FEM with $\gamma=\num{e5}$ and $\lambda=\num{e3}$ for varying melting ranges $\Delta T_f\in\left\{2,1,0.8,0.25,0.175,0.125\right\}\si{K}$. In each case, only about the left half of the enclosure is shown.}
\label{fig:7}
\end{figure}

It is important to note that the melting range $\Delta T_f$ is a non-physical quantity in the sense that gallium, being a pure material with a sharp temperature of fusion, does not posses a mushy region in reality. Therefore, $\Delta T_f$ is a purely numerical value and thus cannot be taken from existing material databases. Instead, different numerical studies are required to deduce an appropriate value for this parameter. In the literature there is no consensus on this choice and we believe the reason for this is that different numerical schemes require a different $\Delta T_f$. However, in the previous subsection it is shown that for the proposed stabilised dG-FEM $\Delta T_f=\SI{0.25}{K}$ yields results in excellent agreement with other research. Based on this situation a parametric study is conducted with the objective of demonstrating the behaviour of the solution when we deviate from this case. Actually, to the authors' knowledge this study is the first attempt to analyse the impact of the melting range $\Delta T_f$ on the flow structure for the problem of gallium melting. \\

In Figure \ref{fig:7} the velocity field and the phase boundary, represented by the $0.5$-contour of $\phi$, at $t=\SI{85}{s}$ computed by the TH type dG-FEM with $\gamma=\num{e5}$ and $\lambda=\num{e3}$ on a $h=\SI{0.4}{mm}$ mesh can be seen for different melting ranges $\Delta T_f$. Obviously, the particular choice of the melting range has a significant impact on the resulting flow structure as both the number and the position of the resulting vortices is affected. We observe that any other choice than $\Delta T_f=\SI{0.25}{K}$ yields solutions which do not possess four separated vortices of approximately equal size, which is understood to be the numerically correct solution. Whilst the upper two vortices are about to merge for $\Delta T_f=\SI{0.8}{K}$, the topmost vortex is smaller than the remaining ones and at the bottom of the enclosure, a fifth vortex weakly appears for $\Delta T_f=\SI{0.175}{K}$. For $\Delta T_f=\SI{1}{K}$ the upper two vortices merge to one, for $\Delta T_f=\SI{2}{K}$ the position of the vortices is wrong and $\Delta T_f=\SI{0.125}{K}$ only yields three vortices. The main reason for the different position and number of vortices is related to the attenuation $A\rb{\phi}$ which acts whenever $\phi<1$. Due to the varying width of the mushy region this attenuation applies differently in each case and yields the apparently different flow structures in the liquid gallium.\\

\begin{figure}[H]
\centering
\subfigure[$\Delta T_f=\SI{2}{K}$]{\includegraphics[width=0.25\textwidth]{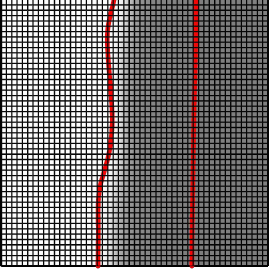}} 		\goodgap	
\subfigure[$\Delta T_f=\SI{1}{K}$]{\includegraphics[width=0.25\textwidth]{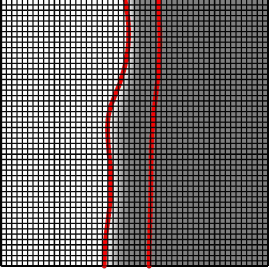}}		\goodgap	
\subfigure[$\Delta T_f=\SI{0.8}{K}$]{\includegraphics[width=0.25\textwidth]{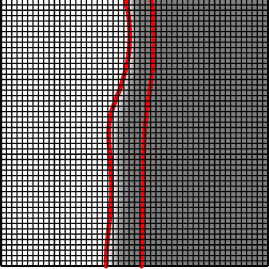}}		\\ 		
\subfigure[$\Delta T_f=\SI{0.25}{K}$]{\includegraphics[width=0.25\textwidth]{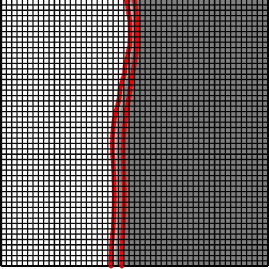}}		\goodgap
\subfigure[$\Delta T_f=\SI{0.175}{K}$]{\includegraphics[width=0.25\textwidth]{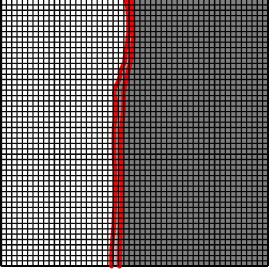}}	\goodgap	 
\subfigure[$\Delta T_f=\SI{0.125}{K}$]{\includegraphics[width=0.25\textwidth]{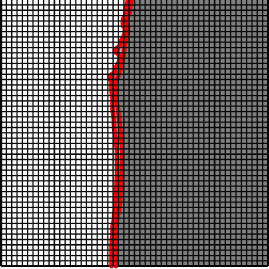}} 	 	
\caption{Underlying $h=\SI{0.4}{mm}$ mesh, phase indicator and mushy region ($0.01$- and $0.99$-contour of $\phi$) at $t=\SI{85}{s}$ computed by the TH type dG-FEM with $\gamma=\num{e5}$ and $\lambda=\num{e3}$ for varying melting ranges $\Delta T_f$. In each case, only a square cutout from the bottom of the enclosure is shown.}
\label{fig:8}
\end{figure}

Lastly, we want to consider the mushy region, being located between the $0.01$- and $0.99$-contour of $\phi$, more closely. Therefore, in Figure \ref{fig:8} a cutout from the bottom of the enclosure together with the underlying $h=\SI{0.4}{mm}$ mesh and the phase indicator function $\phi$ is shown. The light grey represents  the liquid phase while the dark grey indicates the solid phase and the transition between the red lines corresponds to the mushy region. By construction, the width of the mushy region decreases as the melting range $\Delta T_f$ decreases. Note that we decided to present this study only for one fixed mesh size. For this mesh and the TH dG-FEM with $\gamma=\num{e5}$ and $\lambda=\num{e3}$, the smallest value yielding convergent simulations is $\Delta T_f=\SI{0.125}{K}$. Of course, on finer meshes it is possible to simulate the problem for even smaller melting ranges. However, regarding the mushy region for $\Delta T_f=\SI{0.125}{K}$ we observe that the phase transition is already restricted locally to a maximum of two neighbouring elements and is thus located sharply. It is remarkable that the proposed numerical method can deal with the occurring moving interior layers effortlessly even though the mesh is neither adapted to the explicit location of the phase transition, as for example in \cite{DanailaEtAl14}, nor globally refined as for example in \cite{HannounEtAl03}. Summarising, this section shows that stabilised dG-FEM are very well-suited and efficient also for solving thermo-fluid problems involving moving interior layers.

\section{Summary and conclusions}	\label{sec:Conclusions}
In this work we proposed and analysed the performance of a particular class of stabilised dG-FEM for solving thermally-coupled incompressible flow problems with natural convection phenomena based on the Oberbeck--Boussinesq approximation. This class consists of both the mixed-order $\rb{\Quad{-2}\slash\Quad{-1}}\wedge\Quad{-2}$ and the equal-order $\rb{\Quad{-2}\slash\Quad{-2}}\wedge\Quad{-2}$ symmetric interior penalty dG-FEM on quadrilateral meshes with the following two additional stabilisation mechanisms. In order to ensure stability of the equal-order method, pressure jump stabilisation necessarily had to be included whereas the mixed-order method can optionally be equipped with such a term. Most importantly and originally, a classical grad-div stabilisation term has been introduced in combination with all dG-FEM to improve the mass conservation properties of the schemes. The significance of local mass conservation should not be underestimated and is important, for example, in the context of energy balances for indoor airflow simulations. Even though, due to the global full discontinuity of the resulting velocity fields, the proposed class of methods is only approximately $\Hdiv$-conforming, it has been shown that heavy grad-div stabilisation can be used successfully to improve the overall accuracy of the approximate solution in the context of dG-FEM for incompressible natural convection flows. Furthermore, as an additional difficulty we decided to exclusively use non-adapted, uniform meshes to illustrate that the proposed methods are very robust and accurate also for this sub-optimal, but nonetheless very application-relevant situation of having to use under-refined meshes.\\

Therefore, at first the stabilised dG-FEM have been applied to the simulation of the classical differentially heated square cavity for moderate to high Rayleigh numbers as a representative for problems involving both velocity and thermal boundary layers. A detailed qualitative and quantitative analysis by means of comparing Nusselt numbers, stream function values and the fulfilment of the divergence constraint with high-accuracy reference data from the literature has been provided which shows excellent agreement for all considered stabilised dG-FEM. By comparison with a standard conforming FEM and an exactly divergence-free $\Hdiv$-conforming method it turned out that our class of dG-FEM is the superior choice both in terms of accuracy and efficiency. Furthermore, we showed in detail that whilst heavy grad-div stabilisation always improves the mass conservation properties of any finite element type method which is not already exactly divergence-free, for standard conforming FEM the solution deteriorates substantially with an increasing grad-div parameter. The proposed dG-FEM, however, have been shown to not suffer from such a counter-intuitive behaviour. Especially the mixed-order dG-FEM with heavy grad-div and additional pressure jump stabilisation showed the most compelling performance.\\

In the last section we dealt with multiphase flow which classically occurs during melting and solidification processes and involves rather complex moving interior layers. In order to account for such non-isothermal solid/liquid phase transitions the enthalpy-porosity method has been employed. Based on performing best for the heated cavity, for the numerical solution of the resulting mathematical model the grad-div and pressure jump stabilised mixed-order dG-FEM was chosen exemplarily. The problem of melting of pure gallium in a rectangular enclosure has been considered as a benchmark problem for solid/liquid phase change processes. Excellent agreement with previous research has been shown, even though much coarser, non-adapted meshes were used which allow for a more efficient solution of the underlying problem. A mesh sensitivity analysis was provided showing that mesh convergence can be reached relatively fast and that, even on surprisingly coarse meshes, the proposed dG-FEM still yields at least meaningful results. Additionally, a numerical study showed the resulting flow structure for different widths of the melting range from which we deduced that even for remarkably sharp interior layers, the method still converges and remains applicable. \\

Altogether, the proposed class of stabilised dG-FEM performed excellently in all considered studies even though the corresponding problems were highly dynamic, computationally demanding and non-adapted meshes were used. It can thus be inferred that interior penalty dG-FEM, especially in combination with grad-div and pressure jump stabilisation, are highly promising, robust and efficient numerical methods to deal with weakly non-isothermal, natural convection-driven thermo-fluid flows.

\section*{Acknowledgements}	
The authors gratefully acknowledge the helpful comments and suggestions on the manuscript from the anonymous reviewers; they clearly improved the initial version of this work.

\def\bibsection{\section*{References}}

\bibliography{CompPhys_BibTeX}
\bibliographystyle{alphaabbr}

\end{document}